\documentclass[11pt]{article}
\usepackage{amssymb, authblk}
\usepackage{amsmath,bbm}
\usepackage{fullpage}
\usepackage{amsthm, hyperref}
\usepackage{graphicx}
\usepackage{verbatim}
\usepackage{mathtools}
\usepackage{algorithm}
\usepackage{algpseudocode}
\usepackage[dvipsnames]{xcolor}
\usepackage{tikz}
\usetikzlibrary{arrows.meta}
\usetikzlibrary{decorations.pathreplacing}
\DeclareMathAlphabet{\mathpzc}{OT1}{pzc}{m}{it}
\usepackage{natbib}
\usepackage{pdflscape, capt-of}

\algnewcommand\algorithmicinput{\textbf{INPUT:}}
\algnewcommand\INPUT{\item[\algorithmicinput]}
\algnewcommand\algorithmicoutput{\textbf{OUTPUT:}}
\algnewcommand\OUTPUT{\item[\algorithmicoutput]}

\usepackage{hyperref}[]
\hypersetup{
    colorlinks=true,
    linkcolor=blue,
    filecolor=magenta,      
    urlcolor=cyan,
      citecolor=blue,
    }

\usepackage{cleveref}

\newtheorem{theorem}{Theorem}
\newtheorem{lemma}[theorem]{Lemma}

\newtheorem{corollary}[theorem]{Corollary}

\newtheorem{definition}{Definition}

\newtheorem{assumption}{Assumption}
\newtheorem{example}{Example}

\allowdisplaybreaks
\DeclareMathOperator*{\argmin}{arg\,min}
\DeclareMathOperator*{\argmax}{arg\,max}
\DeclareMathOperator*{\sign}{sign}

\title{A review on minimax rates in change point detection and localisation}
\author{Yi Yu}
\affil{Department of Statistics, University of Warwick}
\date{}

\begin{document}

\maketitle

\begin{abstract}

This paper reviews recent developments in fundamental limits and optimal algorithms for change point analysis.  We focus on minimax optimal rates in change point detection and localisation, in both parametric and nonparametric models.  We start with the univariate mean change point analysis problem and review the state-of-the-art results in the literature.  We then move on to more complex data types and investigate general principles behind the optimal procedures that lead to minimax rate-optimal results.
\end{abstract}

\section{Introduction}\label{sec-intro}

Change point analysis, as a statistics research area, can be traced back to the World War II.  \cite{wallis1980statistical} provided a detailed account on how a request from Navy became the prologue of sequential analysis, which can be regarded as a sibling of change point analysis.  The timely demands from manufactory sector during the war boosted the developments of sequential analysis and therefore change point analysis.  The lasting demands from the post-war manufactory sector were a continuing source of fuel fanning the developments of change point analysis in the second half of the 20th century.  In recent years, change point analysis is receiving an unprecedented interest due to the advancement in data collecting, storing and analysing.  We have witnessed the use of change point analysis methods in a wide range of application areas, including finance \citep[e.g.][]{AggarwalInclanLeal1999, AndreouGhysels2002, Ross2013350}, economy \citep[e.g.][]{Fernandez2006}, neuroscience \citep[e.g.][]{Chamroukhi2013633, Lindquist2007CWTDW, RobinsonWL2010}, climatology \citep[e.g.][]{Elsner:2004}, biology \citep[e.g.][]{ErdmanEmerson2008, KwonVannucciSongJeongPfeiffer2008, LioVannucci2000, OliverCarpenaHackenbergBernaola-Galvan2004, PicardLebarbierHoebekeRigaillThiamRobin2011, doi:10.1080/02664763.2013.840272, ShenZhang2012}, chemistry \citep[e.g.][]{ferreira2017partially}, medical sciences \citep[e.g.][]{doi:10.1080/02664763.2013.830085, Henderson/Matthews:1993, doi:10.1080/02664763.2013.809569, McLain_2014}, clinical trials \citep[e.g.][]{KoziolWu1996}, transport \citep[e.g.][]{Hsu1979}, oceanography \citep[e.g.][]{KillickEckleyEwansJonathan2010}, environmental science \citep[e.g.][]{WhitcherByersGuttorpPercival2002}, to name but a few.

Generally speaking, change point analysis is concerned with piecewise-stationary time series data and aims to break them down into stationary pieces.  To be specific, for a length-$T$ time series, we assume that there exists a strictly increasing sequence of unknown time points $\{\eta_k\}_{k = 1}^K \subset \{2, \ldots, T\}$, namely change points, with $K \geq 1$, satisfying that
	\[
		X_t \sim P_t, \quad t \in \mathbb{N}
	\]
	and
	\[
		P_t \neq P_{t-1} \quad \mbox{if and only if } t \in \{\eta_k\}_{k = 1}^K,
	\]
	where $P_t$'s are distributions.  The data $X_t$'s can be scalars, vectors, matrices, networks, functionals, etc.  Given such data, our goal is to estimate the change points accurately.
	
The problem can be further characterised by two additional parameters -- the minimal spacing $\Delta$ and the minimal jump size $\kappa$, which are defined as follows:
	\[
		\Delta = \min_{k = 1, \ldots, K + 1} (\eta_k - \eta_{k-1}) \quad \mbox{and} \quad \kappa = \min_{k = 1, \ldots, K}\kappa_k = \min_{k = 1, \ldots, K} \|P_{\eta_k} - P_{\eta_k - 1}\|_*,
	\]
	where $\|\cdot\|_*$ is a certain distance, $\eta_0 = 1$ and $\eta_{K+1} = T+1$.  Throughout this survey, we let 
	\[
		\kappa \sqrt{\Delta} 
	\]
	be a form of signal-to-noise ratio indicating the fundamental difficulty of the problems.  We remark that this quantity is called differently in different literature, for instance in \cite{verzelen2020optimal} it is called energy.  
	
Given data $\{X_t\}_{t = 1}^T$, we seek estimators $\{\widehat{\eta}_k\}_{k = 1, \ldots, \widehat{K}}$ satisfying that with probability tending to 1 as $T \to \infty$, the following holds:
	\begin{equation}\label{eq-consistent-definition}
		\widehat{K} = K \quad \mbox{and} \quad \lim_{T \to \infty} \frac{\max_{k = 1, \ldots, K} \epsilon_k}{\Delta} = \lim_{T \to \infty} \frac{\max_{k = 1, \ldots, K}|\widehat{\eta}_k - \eta_k|}{\Delta} = 0.
	\end{equation}
	We refer to $\epsilon_k$ the \emph{individual localisation error} of $\eta_k$, $\epsilon = \max_{k = 1, \ldots, K}|\widehat{\eta}_k - \eta_k|$ as the \emph{localisation error} and $\Delta^{-1}\max_{k = 1, \ldots, K}|\widehat{\eta}_k - \eta_k|$ as the \emph{localisation rate}.  For change point estimators satisfying \eqref{eq-consistent-definition}, we call them \emph{consistent} change point estimators. \\

\subsection{What we will cover in this survey}\label{sec-what-we-will-cover}
	
In this survey, we focus on understanding the minimax rates of change point detection and localisation.  These two goals are tightly intertwined.  We endeavour on distinguishing these two concepts in this survey.  Roughly speaking, these two can be regarded as the condition on consistent estimation and the optimal estimation errors.

Rigorously speaking, when $K \geq 1$, the fundamental limit in terms of \textbf{detection} can be presented as a phase transition phenomenon.  One would like to show that in the low signal-to-noise ratio regime
	\begin{equation}\label{eq-intro-detect-lower}
		\kappa \sqrt{\Delta} \lesssim \mbox{ a detection lower bound},
	\end{equation}
	no algorithm is guaranteed to provide consistent change point estimators; and in the high signal-to-noise ratio regime
	\begin{equation}\label{eq-intro-detect-upper}
		\kappa \sqrt{\Delta} \gtrsim \mbox{ a detection upper bound},
	\end{equation}
	we would like to review some computationally-efficient algorithms which can provide consistent change point estimators.
		
On the other hand, the fundamental limit in terms of \textbf{localisation} is that
	\begin{equation}\label{eq-intro-local}
		\inf_{\widehat{\eta}} \sup_{P} \mathbb{E}_P\{d_{\mathrm{H}}(\{\widehat{\eta}_k\}, \{\eta_k\})\} \gtrsim \mbox{optimal localisation error}.
	\end{equation}
	The infimum is taken over all possible estimators of the change points, i.e.~all measurable functions of data.  The supremum is across all possible distributions with signal-to-noise ratios at least higher than optimal localisation conditions.  The distance $d_\mathrm{H}(\cdot, \cdot)$ is the two-sided Hausdorff distance, i.e. for any subset $S_1, S_2 \subset \mathbb{Z}$,
	\[
		d_{\mathrm{H}}(S_1, S_2) = \max\left\{\max_{s_1 \in S_1} \min_{s_2 \in S_2} |s_1 - s_2|, \, \max_{s_2 \in S_2} \min_{s_1 \in S_1} |s_1 - s_2|\right\},
	\]
	with the convention that
	\[
		d_{\mathrm{H}}(S_1, S_2) = \begin{cases}
 			\infty, & S_1 = \emptyset \neq S_2 \mbox{ or } S_2 = \emptyset \neq S_1,\\
 			0, & S_1 = S_2 = \emptyset.
 			\end{cases}
	\]
	The localisation task is to seek change point estimators achieving the optimal localisation rate.

We, ideally, wish for: 
	\begin{itemize}
	\item [(i)] the detection lower and upper bounds in \eqref{eq-intro-detect-lower} and \eqref{eq-intro-detect-upper} coincide in terms of rates; 
	\item [(ii)] the matched detection upper and lower bound coincides with the optimal localisation condition imposed on the distributions considered in \eqref{eq-intro-local}; and 
	\item [(iii)] there exists a computationally-efficient algorithm which can provide estimators with localisation errors matching that in \eqref{eq-intro-local}.  	
	\end{itemize}

In various problems we will cover in this survey, these three goals are achievable simultaneously in some settings, but not all.  We will provide in-depth discussions, including open problems and our conjectures.  In this survey, when talking about optimality, we allow for logarithmic gaps.  When there exist logarithmic gaps, we do not distinguish the term ``nearly-optimal'' and ``optimal''. \\ 

So far, we characterise the changes occurring with an unspecified distance $\|\cdot\|_{*}$ between distinct underlying distributions.  In this survey, we will cover both parametric and nonparametric models.  In parametric models, we will cover univariate mean changes, univariate polynomial coefficients changes, high-dimensional covariance changes, high-dimensional sparse network changes and high-dimensional linear regression coefficient changes.  Absolute values, $\ell_2$-norms of vectors, operator norms and Frobenius norms of matrices, are used as examples of $\|\cdot\|_*$.  In nonparametric models, regarding the distance $\|\cdot\|_{*}$, we will cover the univariate Kolmogorov--Smirnov distance, multivariate supreme norm and a general reproducing kernel Hilbert space distance.  

We will use the univariate mean change problem as the blueprint, studying the fundamental limits of the detection and localisation problems and analysing two types of popular computationally-efficient and statistically-optimal methods.  For all the other aforementioned problems, we will present the information-theoretic lower bounds and an algorithm providing the state-of-art theoretical results.  
		
A summary of these limits can be found in \Cref{tab-summary}.  For detailed parameter definitions, see corresponding sections.  In \Cref{tab-summary}, detection lower bounds correspond to the detection boundaries in the sense of \eqref{eq-intro-detect-lower}, consistence upper bounds correspond to the detection boundaries in the sense of \eqref{eq-intro-detect-upper} and localisation lower bounds correspond to the localisation errors in \eqref{eq-intro-local}.  The optimality upper bounds and localisation upper bounds are the state-of-the art results in the literature.  They are the conditions for computationally-efficient algorithms achieving optimal localisation rates and the localisation errors they achieve.
	
\begin{landscape}
\begin{table}
\caption{Summary.  \label{tab-summary}}
\begin{center}
\renewcommand*{\arraystretch}{3}
	\begin{tabular}{ccccccc}
		Model & $\substack{\mbox{Detection}\\ \mbox{lower bound}}$ & $\substack{\mbox{Consistence}\\ \mbox{upper bound}}$ &  $\substack{\mbox{Optimality}\\ \mbox{upper bound} }$ & $\substack{\mbox{Localisation}\\\mbox{lower bound}}$ & $\substack{\mbox{Localisation}\\ \mbox{upper bound}}$ & Section \\
		\hline
		$\substack{\mbox{Univariate}\\\mbox{mean}}$ & $\substack{\kappa \sqrt{\Delta} \lesssim \\\sigma\log^{1/2}(T)}$ & $\substack{\kappa \sqrt{\Delta} \gtrsim \\\sigma\log^{1/2+\xi}(T)}$ & $\substack{\kappa \sqrt{\Delta} \gtrsim \\\sigma\log^{1/2+\xi}(T)}$ & $\frac{\sigma^2}{\kappa^2}$ & $\frac{\sigma^2 \log(T)}{\kappa^2}$ & \ref{sec-detec-local} \\
		$\substack{\mbox{Univariate} \\ \mbox{polynomials}}$  & $\substack{\kappa \Delta^{r+1/2} \lesssim \\\sigma T^r \log^{1/2}(T)}$ & $\substack{\kappa \Delta^{r+1/2} \gtrsim \\ \sqrt{K}\sigma T^r \log^{1/2+\xi}(T)}$ & $\substack{\kappa \Delta^{r+1/2} \gtrsim \\\sqrt{K} \sigma T^r \log^{1/2+\xi}(T)}$ & $\left(\frac{\sigma^2}{\kappa^2 T}\right)^{\frac{1}{2r+1}}$ & $\left\{\frac{\sigma^2\log(T)}{\kappa^2T}\right\}^{\frac{1}{2r+1}}$ & \ref{sec-piecewise-polynomial} \\
		$\substack{\mbox{High-dim} \\ \mbox{covariance}}$ & $\substack{\kappa \sqrt{\Delta} \lesssim \\ \sigma^2\sqrt{p}}$ & $\substack{\kappa \sqrt{\Delta} \gtrsim \\\sigma^2 \sqrt{p}\log^{1/2+\xi}(T)}$ & $\substack{\kappa \sqrt{\Delta} \gtrsim \\\sigma^2 \sqrt{p}\log^{1/2+\xi}(T)}$ & $\frac{\sigma^4}{\kappa^2}$ &  $\frac{\sigma^4 \log(T)}{\kappa^2}$ & \ref{sec-cov-change} \\
		$\substack{\mbox{Dynamic} \\ \mbox{networks}}$ & $\substack{\kappa_0 \sqrt{\Delta} \lesssim \\ \log^{1/2}(T)/\sqrt{n\rho}}$ & $\substack{\kappa_0 \sqrt{\Delta} \gtrsim \\ \log^{1 + \xi}(T)/\sqrt{n\rho}}$ & $\substack{\kappa_0 \sqrt{\Delta} \gtrsim \\ \sqrt{r}\log^{1 + \xi}(T)/\sqrt{n\rho}}$ & $\frac{1}{\kappa_0^2n^2\rho}$ & $\frac{\log^2(T)}{\kappa_0^2n^2\rho}$ & \ref{sec-graphon}\\
		$\substack{\mbox{High-dim} \\ \mbox{linear} \\ \mbox{regression}}$ & $\substack{\kappa \sqrt{\Delta} \lesssim \\ \sigma_{\varepsilon}\sqrt{d_0 \log^{1+\xi}(T)}}$ & $\substack{\kappa \sqrt{\Delta} \gtrsim \\ d_0 \sigma_{\varepsilon}\sqrt{K \log^{1+\xi}(T)}}$ & $\substack{\kappa \sqrt{\Delta} \gtrsim \\ d_0 \sigma_{\varepsilon}\sqrt{K \log^{1+\xi}(T)}}$ & $\frac{d_0\sigma_{\varepsilon}^2}{\kappa^2}$ & $\frac{d_0\sigma_{\varepsilon}^2\log(T \vee p)}{\kappa^2}$ & \ref{sec-regression} \\		
		$\substack{\mbox{Univariate} \\ \mbox{nonparametric}}$ & $\kappa\sqrt{\Delta} \lesssim 1$ & $\substack{\kappa\sqrt{\Delta} \gtrsim \\ \log^{1/2+\xi}(T)}$ & $\substack{\kappa\sqrt{\Delta} \gtrsim \\ \log^{1/2+\xi}(T)}$ & $\kappa^{-2}$ & $\frac{\log(T)}{\kappa^2}$ & \ref{sec-np-uni} \\		
		$\substack{\mbox{Multivariate} \\ \mbox{nonparametric}}$ & $\kappa^{p/2+1}\sqrt{\Delta} \lesssim 1$ & $\substack{\kappa^{p/2+1}\sqrt{\Delta} \gtrsim \\\log^{1/2+\xi}(T)}$ &  $\substack{\kappa^{p/2+1}\sqrt{\Delta} \gtrsim \\\log^{1/2+\xi}(T)}$ & $\kappa^{-(p+2)}$ & $\frac{\log(T)}{\kappa^{p+2}}$ & \ref{sec-np-multi}\\
		$\substack{\mbox{An RKHS} \\ \mbox{example}}$ & $\substack{\kappa \sqrt{\Delta} \lesssim \\ \log^{1/2}(T)}$ & $\substack{\kappa \sqrt{\Delta} \gtrsim \\ K \log^{1/2+\xi}(T)}$ & $\substack{\kappa \sqrt{\Delta} \gtrsim \\ K \log^{1/2+\xi}(T)}$ & $\kappa^{-2}$ & $\frac{\log(T)}{\kappa^2}$ & \ref{sec-general-np} \\
	\end{tabular}	
\end{center}
\end{table}
\end{landscape}

\subsection{What we will not cover in this survey}\label{sec-what-not-cover}	

After more than seven decades of developing, change point analysis has become an extremely fruitful area with numerous papers from a wide range of areas.  Even within the statistics community, change point analysis has been studied extensively from many different angles.  There are a few important topics we do not intend to discuss in depth in this survey.  We briefly mention them here.

\medskip
\textbf{Online change point analysis.}  In this survey, we will only focus on offline change point analysis, i.e.~given data $\{X_t\}_{t = 1}^T$, we retrospectively seek change points in $\{1, \ldots, T\}$.  Another important area of change point analysis is called online/sequential change point analysis, in which one is making sequential decisions on whether a change point has occurred while collecting data.  

The statistical problems associated with online change point detection include minimising the detection delay, e.g.~upper bounding $(\widehat{\eta} - \eta)_+$, while controlling false positives, e.g.~upper bounding the probability of $\widehat{\eta} < \eta$.  There is a vast body of existing literature on this topic, including \cite{moustakides1986optimal}, \cite{ritov1990decision}, \cite{lorden1971procedures}, \cite{lai1981asymptotic}, \cite{lai1995sequential}, \cite{lai1998information}, \cite{lai2001sequential}, \cite{lai2010sequential}, \cite{chu1996monitoring}, \cite{aue2004delay}, \cite{aue2009delay}, \cite{kirch2008bootstrapping}, \cite{huvskova2012bootstrapping}, \cite{mei2010efficient}, \cite{huvskova2010fourier}, \cite{hlavka2016bootstrap}, \cite{desobry2005online}, \cite{fearnhead2007line}, \cite{he2018sequential}, \cite{kirch2018modified}, \cite{kurt2018real}, \cite{chen2019sequential}, \cite{dette2019likelihood}, \cite{gosmann2019new}, \cite{dette2019likelihood}, \cite{keshavarz2018sequential}, \cite{chen2020high}, \cite{siegmund2013sequential},  \cite{tartakovsky2014sequential}, \cite{namoano2019online}, \cite{maillard2019sequential} and \cite{yu2020note}, among others.

\medskip
\textbf{Testing.}  Testing and estimation are two indispensable pillars in statistical problems.  In \Cref{sec-what-we-will-cover}, we mentioned that detection and localisation are two different estimation tasks in change point analysis.  The testing aspects in change point analysis focus on the Type-I and -II errors controls on testing the presence of change points, and also include the limiting distributions of change point estimators, constructing confidence intervals of change points, etc.  As a statistical problem, testing is generally easier than estimation, in terms of the fundamental limits.  Quite the contrary, the study of the fundamental limits of the testing problem is lagged behind.  The literature on different aspects of testing includes \cite{yao1989least}, \cite{frick2014multiscale}, \cite{enikeeva2019bump}, \cite{vanegas2019multiscale}, \cite{dette2019detecting}, \cite{dette2018multiscale}, \cite{akashi2018change}, \cite{dette2018change}, \cite{aue2018detecting}, \cite{aue2013structural}, \cite{robbins2011mean}, \cite{liu2019minimax}, \cite{stoehr2020detecting}, \cite{kirch2015detection}, \cite{jewell2019testing}, \cite{chen2019change}, \cite{Jirak2015}, \cite{chu2019asymptotic} and \cite{verzelen2020optimal}, among others. 

\medskip
\textbf{Computation.}  In this survey, for each problem, we will only present one or two polynomial-time algorithms which provide nearly-optimal results.  In practice, especially for high-dimensional data, it is crucial to improve the computational efficiency without sacrificing too much statistical accuracy.  There is a line of attack on improving the computational time of the methods we will introduce later in this survey.  These works include \cite{romano2020detecting}, \cite{hocking2020constrained}, \cite{tickle2020parallelization}, \cite{haynes2017computationally}, \cite{hocking2017log}, \cite{maidstone2017optimal}, \cite{haynes2017computationally-2}, \cite{killick2012optimal}, \cite{rigaill2010pruned}, \cite{kovacs2020seeded} and \cite{kovacs2020optimistic}, among others.

\medskip
\textbf{Tuning parameter selection.}  For all the methods studied in this survey, the theoretical results rely on some properly chosen tuning parameters.  This is always an important but hard-to-address problem in statistics.  Generally speaking, some papers use information-type criteria pioneered in \cite{yao1988estimating}, and others use data-driven methods \citep[e.g.][]{padilla2019optimal, matteson2014nonparametric}.  It is fair to say that tuning parameter selection is still an open topic in change point analysis.

\section{Univariate mean changes: a blueprint}\label{sec-detec-local}

\subsection{Setup and overview}

Arguably, the simplest and best-studied change point analysis problem is the univariate time series with piecewise-constant mean and independent sub-Gaussian noise.  We formalise the problem below.

\begin{assumption}\label{assume:change-uni-mean}
Let $\{X_t\}_{t = 1}^T \subset \mathbb{R}$ be independent sub-Gaussian random variables with continuous density such that $\mathbb{E}(X_t) = f_t$ and $\|X_t\|_{\psi_2} \leq \sigma$ for all $t \in \{1, \ldots, T\}$.  

Let $\{\eta_k\}_{k=0}^{K+1} \subset \{1, \ldots, T+1\}$ be a collection of change points such that $1 = \eta_0 < \eta_1 < \ldots < \eta_K \leq T < \eta_{K+1} = T+1$ and $f_t \neq f_{t - 1}$, if and only if $t \in \{\eta_k\}_{k = 1}^K$.

Assume the minimal spacing $\Delta$ and the jump size $\kappa$ are defined to be
	\[
		\Delta = \min_{k = 1, \ldots, K+1} \bigl\{\eta_k - \eta_{k-1}\bigr\} > 0,
	\]	
	and
	\[
		\kappa = \min_{k = 1, \ldots, K} \kappa_k = \min_{k = 1, \ldots, K} \bigl|f_{\eta_k} - f_{\eta_k-1}\bigr|> 0.
	\]
\end{assumption} 

Note that the $\|\cdot\|_{\psi_2}$ is the Orlicz-$\psi_2$-norm or the sub-Gaussian norm, defined as, for any random variable $X$,
	\[
		\|X\|_{\psi_2} = \inf\left\{t > 0: \, \mathbb{E}\left\{\exp(X^2/t^2)\right\} \leq 2\right\}.
	\]
	The condition on the continuous density is merely to impose uniqueness of the estimators.  We impose it here just for simplicity.
	
As for the problem detailed in \Cref{assume:change-uni-mean}, the detection lower bound is $\sigma \log^{1/2}(T)$ and the localisation lower bound is $\sigma^2 \kappa^{-2}$.  These two results are in Lemmas~\ref{lemma-low-snr} and \ref{lemma-error-opt}, respectively.  To match these lower bounds, we will show two nearly-optimal polynomial-time algorithms. 

\subsection{Detection boundary}\label{sec-uni-detect}
		
The detection boundary $\sigma \log^{1/2}(T)$ has been established in various different papers, including \cite{chan2013}, \cite{frick2014multiscale}, \cite{dumbgen2001multiscale}, \cite{dumbgen2008multiscale}, \cite{LiEtal2017}, \cite{jeng2012simultaneous}, \cite{enikeeva2018bump} and \cite{wang2020univariate}.  We formalise the result below.

\begin{lemma}[Lemma 1 in \citealp{wang2020univariate}]\label{lemma-low-snr}
Let $\{X_t\}_{t=1}^T$ be a time series satisfying \Cref{assume:change-uni-mean}.  Let $P^T_{\kappa, \Delta, \sigma}$ denote the corresponding joint distribution.  For any $0 < c < 1$, consider the class of distributions
	\[
		\mathcal{P}^T = \left\{P^T_{\kappa, \Delta, \sigma}: \, \Delta = \min \left\{\left\lfloor c\frac{\log(T)}{\kappa^2/\sigma^2}\right\rfloor, \left\lfloor \frac{T}{4} \right\rfloor\right\} \right\}.
	\]
	Then, there exists an $T(c)$, which depends on $c$, such that, for all $T$ larger than $T(c)$, 
	\[
		\inf_{\{\widehat{\eta}\}} \sup_{P \in \mathcal{P}^T} \mathbb{E}_P\bigl\{d_{\mathrm{H}}(\{\widehat{\eta}\}, \, \{\eta(P)\} ) \bigr\} \geq \frac{\Delta}{2},
	\]
	where the infimum is over all estimators $\{\widehat{\eta}\}$ of the change point locations and $\{\eta(P)\}$ is the set of locations of the change points of $P \in \mathcal{P}^T$.
\end{lemma}

For the localisation rate achieved in \Cref{lemma-low-snr}, it holds that 
	\[
		\frac{d_{\mathrm{H}}(\{\widehat{\eta}\}, \{\eta\})}{\Delta} \geq \frac{1}{2}, 
	\]
	which does not vanish.  Corresponding to \eqref{eq-consistent-definition}, \Cref{lemma-low-snr} shows that if the signal-to-noise ratio is in the regime
	\[
		\kappa\sqrt{\Delta} \lesssim \sigma \log^{1/2}(T),
	\]
	then no algorithm is guaranteed to provide consistent estimators.  

To complete the phase transition phenomenon, one needs to provide computationally-efficient algorithms, which provide consistent change point estimators in the regime
	\[
		\kappa\sqrt{\Delta} \gtrsim \sigma \log^{1/2}(T).
	\]

As for change point analysis, there are two main types of algorithms, which we will refer to as penalisation-based estimators and scan-statistics-type estimators.  The penalisation-based estimators are based on a penalised objective function.  The objective function is often a goodness-of-fit criterion and the penalty is usually imposed on the number of change points, to avoid overfitting.  The scan-statistics-type estimators are obtained by adopting a certain statistic, scanning through most if not all the data to evaluate each time point's potential of being a change point.  Each of these two types of estimation contains a variety of methods, and each of these two has a representative method being nearly-optimal in the sense we focus on in this paper.  
	
\subsubsection{Penalisation-based estimators}\label{sec-uni-mean-ell-0}

Let $\mathcal{P}$ be any {\it interval partition} of $\{1, \ldots, T\}$, i.e.~a collection of $|\mathcal{P}| \geq 1$ disjoint interval subsets of $\{1, \ldots, T\}$ in the form of
	\[
		\mathcal{P} = \bigl\{\{1, \ldots, i_1 - 1\}, \{i_1, \ldots, i_2-1\}, \ldots, \{i_{|\mathcal{P}|-1}, \ldots, i_{|\mathcal{P}|}-1 \}\bigr\},
	\]
	for some integers $1 < i_1 < \cdots < i_{|\mathcal{P}|-1} \leq T < i_{|\mathcal{P}|} = T+1$.  For a positive tuning parameter $\lambda > 0$ and data $\{X_t\}_{t=1}^T$, let 
	\begin{equation}\label{eq-p-hat}
		\widehat{\mathcal{P}}(\lambda) \in \argmin_{\mathcal P} G\bigl(\mathcal P, \{X_t\}_{t=1}^T, \lambda\bigr),
	\end{equation}
	where the minimum ranges over all interval partitions of $\{1, \ldots, T\}$ and, for any such partition $\mathcal{P}$,
	\begin{equation}\label{eq-G}
		G\bigl(\mathcal P, \{X_t\}_{t=1}^T, \lambda\bigr) = \sum_{I\in \mathcal P} H(I) + \lambda |\mathcal{P}| =  \sum_{I\in \mathcal P} \|X_I - \Pi_I X_I\|^2 + \lambda |\mathcal{P}|,
	\end{equation}
	where $\Pi_I$ is the projection matrix of the subspace spanned by an all-one vector, i.e.
	\begin{equation}\label{eq-pi-i-uni-mean}
		\Pi_I = \mathbbm{1}_I \left(\mathbbm{1}_I^{\top} \mathbbm{1}_I\right)^{-1} \mathbbm{1}_I^{\top}, \quad \mathbbm{1}_I = (1, \ldots, 1)^{\top} \in \mathbb{R}^{|I|}
	\end{equation}
	and $X_I = (X_i, i \in I)^{\top} \in \mathbb{R}^{|I|}$.  In fact, for any interval $I \subset \{1, \ldots, T\}$, 
	\[
		\Pi_I X_I = |I|^{-1} \sum_{i \in I} X_i.
	\]
	We adopt the seemingly unnecessary notation \eqref{eq-pi-i-uni-mean} to be consistent with that in \Cref{sec-piecewise-polynomial}.

The optimization problem \eqref{eq-p-hat} is known as the minimal partition problem and can be solved using dynamic programming in polynomial time \citep[e.g. Algorithm~1 in][]{FriedrichEtal2008}.  The change point estimator resulting from the solution to \eqref{eq-p-hat} is simply obtained from taking all the left endpoints of the intervals  $I \in \widehat{\mathcal{P}}$, while letting $\widehat{\eta}_0 = 1$.  For completeness, we include the algorithm in \Cref{algorithm:PDP}.  As we have emphasised in \Cref{sec-what-not-cover}, the computational issue is not covered in this paper.  The computational cost of \Cref{algorithm:PDP} is of order $O(T^2 \mathrm{Cost}(T))$, where $\mathrm{Cost}(T)$ is the computational cost of calculating the function $H(I)$, with an interval $I$ of length $T$.  There exist more efficient variants of \Cref{algorithm:PDP} in solving the optimisation problem \eqref{eq-p-hat}, including \cite{killick2012optimal}.

\begin{algorithm}[htbp]
\begin{algorithmic}
	\INPUT Data $\{X(t)\}_{t=1}^{T}$, tuning parameter $\lambda > 0$.
	\State $(\mathcal{B}, s, t, \mathrm{FLAG}) \leftarrow (\emptyset, 0, 2, 0)$
	\While{$s < T-3$}
		\State $s \leftarrow s + 1$
		\While{$t < T$ and $\mathrm{FLAG} = 0$}
			\State $t \leftarrow t+1$ 			
			\If{$\min_{l \in \{s+1, \ldots, t-1\}}\left\{H([s, l]) + H([l+1, t])\right\} + \lambda <  H([s, t])$} 
				\State $s \leftarrow \min\Big\{l \in s+1, \ldots, t-1: \, H([s, l]) + H([l+1, t]) + \lambda <  H([s, t]) \Big\}$
				\State $\mathcal{B} \leftarrow \mathcal{B} \cup \{s\}$
				\State $\mathrm{FLAG} \leftarrow 1$
			\EndIf
		\EndWhile
	\EndWhile
	\OUTPUT The set of estimated change points $\mathcal{B}$.
\caption{Penalised dynamic programming. }
\label{algorithm:PDP}
\end{algorithmic}
\end{algorithm}

\begin{theorem}[Theorem 3 in \citealp{wang2020univariate}]\label{prop:1d localization}
Let $\{X_i\}_{i=1}^T$ satisfy \Cref{assume:change-uni-mean}.  Assume that there exists a sufficiently large absolute constant $C_{\mathrm{SNR}} > 0$ such that for any $\xi > 0$,
	\[
		\kappa\sqrt{\Delta}/\sigma \geq C_{\mathrm{SNR}}\sqrt{\log^{1+\xi}(T)}.
	\]
	For any $\lambda>0$, let $\{\widehat{\eta}_k\}_{k = 1, \ldots, \widehat{K}}$ be the output of \Cref{algorithm:PDP} with function $H(\cdot)$ defined in \eqref{eq-G}.  We have that, for any choice of $c>0$, there exists a constant  $C_{\lambda} > 0$, which depends on $c$ such that, for $\lambda = C_{\lambda}\sigma^2 \log(T)$, it holds that
	\[
		\mathbb{P} \bigl\{\widehat{K} = K \quad \mbox{and} \quad \epsilon_k = |\widehat \eta_k(\lambda) -\eta_k| \leq C_{\epsilon}\sigma^2 \log(T) /\kappa^2_k, \, \forall k \in \{1,\ldots,K\} \bigr\} \geq 1 - T^{-c},
	\]
	where $C_{\epsilon} > 0$ is a constant depending on $C_{\lambda}$ and $C_{\mathrm{SNR}}$.  
\end{theorem}

\Cref{prop:1d localization} shows that in the regime $\kappa\sqrt{\Delta} \gtrsim \sigma \log^{1/2+\xi}(T)$, the outputs of \Cref{algorithm:PDP} with \eqref{eq-G} are consistent.  To be specific, in this signal-to-noise ratio regime,
	\[
		\lim_{T \to \infty}\frac{\max_{k = 1, \ldots, K}\epsilon_k}{\Delta} \lesssim \lim_{T \to \infty}\frac{\max_{k = 1, \ldots, K}\sigma^2 \log(T) /\kappa^2_k}{\Delta} \lesssim \lim_{T\to \infty} \log^{-\xi}(T) = 0.
	\]
	This also explains the role of $\xi$.  It is introduced merely for mathematical purposes on enforcing the vanishing ratio and for notational simplicity.  The term $\log^{\xi}(T)$ can be replaced by any diverging sequence $a_T$.

\Cref{prop:1d localization} and \Cref{lemma-low-snr} together show a phase transition phenomenon that:
	\begin{itemize}
	\item in the low signal-to-noise ratio regime
		\[
			\kappa\sqrt{\Delta} \lesssim \sigma \log^{1/2}(T),
		\]	
		no algorithm is guaranteed to be consistent; and
	\item in the high signal-to-noise ratio	regime
		\[
			\kappa\sqrt{\Delta} \lesssim \sigma \log^{1/2 + \xi}(T), \quad \forall \xi > 0,
		\]
		we have a computationally-efficient algorithm which achieves a consistent change point estimation.
	\end{itemize}

We remark that there are other types of penalisations.  In \eqref{eq-G}, the penalty is imposed on the number of change points and is equivalent to an $\ell_0$ penalty.  It is natural to replace the $\ell_0$ penalty with an $\ell_1$ penalty and ends up with a fused Lasso \citep{tibshirani2005sparsity} or a trend filtering \citep[e.g.][]{tibshirani2014adaptive} problem.  There have indeed been works analysing change points using $\ell_1$ penalties due to its computational efficiency.  It is known that in terms of change point detection and localisation, $\ell_1$ penalisation based methods are sub-optimal  \citep{lin2016approximate}, but the estimators can be improved with proper post-processing \citep[e.g.][]{zhang2019element, hyun2018exact}.

\subsubsection{Scan-statistics-type estimators}	

Arguably, the most popular statistic used in change point analysis is the cumulative sum \citep[CUSUM,][]{Page1954} statistic, which was proposed as an extension of the sequential probability ratio test statistics \citep{Wald1945}.  

\begin{definition}[CUSUM statistics]\label{def-cusum}
For a sequence $\{X_t\}_{t=1}^T$, any integer triplet $(s, t, e)$, $0 \leq s < t < e \leq T$, let the CUSUM statistic be
	\[
		\widetilde X^{s, e}_t = \sqrt {\frac{e-t}{(e-s)(t-s)} } \sum_{i=s+1}^t X_i  -\sqrt {\frac{t-s}{(e-s)(e-t)} } \sum_{i=t+1}^e X_i.
	\]
\end{definition}

We will encounter multiple versions of \Cref{def-cusum} in the rest of this survey.  The original CUSUM statistic is restricted to the case that $X_t$'s are scalars, but they will be allowed to be in different spaces in this survey.  \\

The CUSUM statistic is originated from a log-likelihood ratio test statistic.  We elaborate this from the example below.

\begin{example}\label{ex-1}
Let $\{X_t\}_{t = 1}^T$ be a sequence of independent Gaussian random variables with unknown mean $\mu_t$ and known variance $\sigma^2$.  For any $t \in \{1, \ldots, T-1\}$, let
	\begin{align*}
		H_{1, t}: \mu_1 = \cdots = \mu_t \neq \mu_{t+1} = \cdots = \mu_T.
	\end{align*}
	We want to test
	\begin{align}\label{eq-cusum-test-two-sample}
		H_0: \mu_1 = \cdots = \mu_T \quad \mbox{vs.} \quad H_1 = \cup_{t = 1}^{T-1} H_{1, t}.
	\end{align}
\end{example}

For any fixed $t \in \{1, \ldots, T-1\}$, define
	\[
		\overline{X}_1 = \frac{1}{t}\sum_{i = 1}^t X_i, \quad \overline{X}_2 = \frac{1}{T-t}\sum_{i = t+1}^T X_i \quad \mbox{and} \quad \overline{X} = \frac{1}{T}\sum_{i = 1}^T X_i.
	\]
	The generalised likelihood ratio test statistic of the problem in \Cref{ex-1} is that 
	\begin{align*}
		T_t & = \log\left[\frac{\prod_{i=1}^t \frac{1}{\sqrt{2\pi \sigma^2}} \exp\left\{-\frac{(X_i - \overline{X}_1)^2}{2\sigma^2}\right\} \prod_{i=t+1}^T \frac{1}{\sqrt{2\pi \sigma^2}} \exp\left\{-\frac{(X_i - \overline{X}_2)^2}{2\sigma^2}\right\} }{\prod_{i=1}^T \frac{1}{\sqrt{2\pi \sigma^2}} \exp\left\{-\frac{(X_i - \overline{X})^2}{2\sigma^2}\right\}} \right]	 \\
		& = \frac{1}{2\sigma^2} \left\{\sum_{i = 1}^T (X_i - \overline{X})^2 - \sum_{i = 1}^t (X_i - \overline{X}_1)^2 - \sum_{i = t + 1}^T (X_i - \overline{X}_2)^2\right\} \\
		& = \frac{1}{2\sigma^2} \left\{\sqrt{\frac{T-t}{Tt}} \sum_{i = 1}^t X_i - \sqrt{\frac{t}{T(T-t)}} \sum_{i = t+1}^T X_i\right\}^2 = \frac{1}{2\sigma^2} \left(\widetilde{X}_t^{0, T}\right)^2.
	\end{align*}
Then \eqref{eq-cusum-test-two-sample} can be conducted based on 
	\[
		\max_{t = 1, \ldots, T-1} T_t = \frac{1}{2\sigma^2}\max_{t = 1, \ldots, T-1} \left(\widetilde{X}_t^{0, T}\right)^2,
	\]
	which is equivalent to the use of CUSUM statistics in change point detection.  
	
In fact, CUSUM statistics can be used for change point detection in a number of ways.  Arguably, the most popular and standard method is the binary segmentation \citep[e.g.][]{ScottKnott1974, venkatraman1992consistency, vostrikova1981detection}.  The key idea is to find 
	\[
		\widehat{t} \in \max_{t = 1, \ldots, T-1} |\widetilde{X}^{0, T}_t|.
	\]
	For a pre-specified threshold $\tau$, if $|\widetilde{X}^{0, T}_{\widehat{t}}| \geq \tau$, then we declare $\widehat{t}$ to be a change point estimator and the procedure is conducted on the intervals $(0, \widehat{t})$ and $[\widehat{t}, T]$ respectively.  The procedure is terminated if there is no more change point estimator declared, or if the resulting interval is too narrow.  Binary segmentation is a computationally-efficient algorithm, but sub-optimal.  The sub-optimality can be intuitively explained as follows.  When there are potentially multiple change points, the consecutive change points may cancel out each other \citep[e.g.][]{fryzlewicz2014wild}.
	
In order to improve the theoretical guarantees of CUSUM-based algorithms, especially to tackle the multiple change points scenario, a large number of variants have been proposed, including \cite{fryzlewicz2014wild}, \cite{kovacs2020seeded}, \cite{kovacs2020optimistic}, \cite{anastasiou2019detecting}, \cite{baranowski2016narrowest}, among others.

We use the wild binary segmentation \citep{fryzlewicz2014wild} as an example to illustrate how a CUSUM-based method can achieve optimality. 

\begin{algorithm}[htbp]
\begin{algorithmic}
	\INPUT Independent samples $\{X_t\}_{t=1}^T$, collection of intervals $\{(\alpha_m, \beta_m)\}_{m=1}^M$, tuning parameter~$\tau > 0$.
	\For{$m = 1, \ldots, M$}  
		\State $(s_m, e_m) \leftarrow (s, e)\cap [\alpha_m, \beta_m]$
		\If{$e_m - s_m > 1$}
			\State $b_{m} \leftarrow \argmax_{t = s_m + 1, \ldots, e_m - 1}| \widetilde{X}^{s_m, e_m}_t|$
			\State $a_m \leftarrow \bigl|\widetilde{X}^{s_m, e_m}_{b_m}\bigr|$
		\Else 
			\State $a_m \leftarrow -1$	
		\EndIf
	\EndFor
	\State $m^* \leftarrow \argmax_{m = 1, \ldots, M} a_m$
	\If{$a_{m^*} > \tau$}
		\State add $b_{m^*}$ to the set of estimated change points
		\State WBS$((s, b_{m*}),\{ (\alpha_m,\beta_m)\}_{m=1}^M, \tau)$
		\State WBS$((b_{m*}+1,e),\{ (\alpha_m,\beta_m)\}_{m=1}^M,\tau ) $

	\EndIf  
	\OUTPUT The set of estimated change points.
\caption{Wild Binary Segmentation. WBS$((s, e),$ $\{ (\alpha_m,\beta_m)\}_{m=1}^M, \tau $)}
\label{algorithm:WBS}
\end{algorithmic}
\end{algorithm}

\begin{theorem}[Theorem 4 in \citealp{wang2020univariate}]\label{thm-wbs-uni}
Assume that the inputs of \Cref{algorithm:WBS} are as follows.
	\begin{itemize}
		\item The sequence $\{X_t\}_{t=1}^T$ satisfies \Cref{assume:change-uni-mean}.  In addition, assume that there exists a sufficiently large absolute constant $C_{\mathrm{SNR}} > 0$ such that for any $\xi > 0$,
		\[
			\kappa\sqrt{\Delta}/\sigma \geq C_{\mathrm{SNR}}\sqrt{\log^{1+\xi}(T)}.
		\]
		\item The collection of intervals $\{(\alpha_m, \beta_m)\}_{m=1}^M \subset \{1, \ldots, T\}$, whose endpoints are drawn independently and uniformly from $\{1, \ldots, T\}$, satisfy 
		\[
		\max_{m = 1, \ldots, M}(\beta_m - \alpha_m) \leq C_R \Delta, 
		\]
		almost surely, for an absolute constant $C_{R} > 1$.
		\item The tuning parameters  $\tau$ satisfies 
			\[
				c_{\tau, 1}\sigma\sqrt{\log(T)} < \tau < c_{\tau, 2}\kappa\sqrt{\Delta},
			\]
			where $c_{\tau, 1}, c_{\tau, 2} > 0$ are sufficiently large and   small  absolute constants.
	\end{itemize}	
	Let $\bigl\{\widehat{\eta}_k\bigr\}_{k=1}^{\widehat{K}}$ be the corresponding output of \Cref{algorithm:WBS}.  It holds that 
	\begin{align}
	& \mathbb{P}\left\{\widehat{K} = K \quad \text{and} \quad   \epsilon_k = |\widehat{\eta}_k - \eta_k| \leq C_{\epsilon}\sigma^2\log(T)\kappa^{-2}_k, \forall k \in \{1,\ldots,K\} \right\} \nonumber \\
	& \hspace{2cm} \geq 1 - T^{-c} - \exp\left\{\log\left(\frac{T}{\Delta}\right) - \frac{M\Delta^2}{16T^2}\right\}, \nonumber 
	\end{align}
	where $C_{\epsilon}, c > 0$ are absolute constants.
\end{theorem}

Following the same discussions after \Cref{prop:1d localization}, provided $\log(T/\Delta) \lesssim M\Delta^2 T^{-2}$, \Cref{thm-wbs-uni} shows that \Cref{algorithm:WBS} provides consistent change point estimation under a nearly-optimal signal-to-noise ratio regime.  The key to the success of \Cref{algorithm:WBS} is the usage of random intervals, but in order to achieve the optimality, in \Cref{thm-wbs-uni}, the lengths of the random intervals are at most of the order of the minimal spacing.  This is of course not practical, but essential in deriving the optimality.  Similar treatments can be found in other forms, such as the parameter $\beta$ used in \cite{wang2016high}.  If we relax the condition that $C_R$ being an absolute constant, then $C_R \leq T/\Delta$.  This results in an inflation in the required signal-to-noise ratio and the resulting localisation rate.  To be specific, one would require
	\[
		\kappa\sqrt{\Delta}/\sigma \geq C_{\mathrm{SNR}}\frac{T}{\Delta}\sqrt{\log^{1+\xi}(T)}
	\]
	and have the localisation rate being 
	\[
		C_{\epsilon}\sigma^2\log(T)\kappa^{-2}_k\frac{T^2}{\Delta^2}.
	\]

As we have pointed out, the sub-optimality of the binary segmentation roots in the multiple change point scenario.  Note that the detection upper bound is $\kappa \sqrt{\Delta} \asymp \sigma  \log^{1/2+\xi}(T)$.  For simplicity, we let $\kappa, \sigma \asymp 1$ and $\xi = 1/2$, then this means $\Delta$ can be as small as $\log^2(T)$.  In this case, the number of change points can be as many as $T/\Delta \asymp T \log^{-2}(T)$, which diverges as $T$ grows unbounded.  

Most if not all of the variants of the binary segmentation works on how to narrow the focus to intervals containing only finite number of true change points.  Different variants use different additional parameters to guarantee this for theoretical purposes.  To the best of our knowledge, there is no algorithm can deal with this issue satisfactorily both theoretically and practically.  For example, the WBS-type methods require this additional constant $C_R$ in the upper bound on the lengths of random intervals.  The narrowest-over-threshold method \citep{baranowski2016narrowest} is shown to be too sensitive to tuning parameters in numerical experiments.  The optimistic search strategy \citep{kovacs2020optimistic} works under a stronger condition on the minimal spacing for the multiple change points scenario.

The CUSUM statistics essentially can be regarded as differences between weighted sample means.  The weights play the role of variance stabilisation.  We remark that there are other types of scan statistics, including those used in \cite{CribbenYu2017}, \cite{LiuEtal2018} and \cite{niu2012screening}.  

\subsection{Optimal localisation rate}\label{sec-uni-local}

As for the localisation, we have the following minimax lower bound.

\begin{lemma}[Lemma 2 in \citealp{wang2020univariate}]\label{lemma-error-opt}
Assume that the sequence $\{X_t\}_{t=1}^T$ satisfies \Cref{assume:change-uni-mean}.  Let $P^T_{\kappa, \Delta, \sigma}$ denote the corresponding joint distribution.  Consider the class of distributions
	\[
		\mathcal{Q}^T = \left\{P^T_{\kappa, \Delta, \sigma}: \, \Delta < T/2, \, \kappa\sqrt{\Delta}/\sigma \geq \zeta_T \right\},
	\]
	for any sequence $\{\zeta_T\}$ such that $\lim_{T \rightarrow \infty} \zeta_T = \infty $.  Then, for all $T$ large enough, it holds that 
	\[
		\inf_{\hat{\eta}} \sup_{P \in \mathcal{Q}^T} \mathbb{E}_P\bigl(\bigl|\widehat{\eta} - \eta(P)\bigr|\bigr) \geq \max \left\{ 1, \frac{1}{2} \Big\lceil\frac{\sigma^2}{\kappa^2} \Big\rceil e^{-2} \right\}, 
	\]
	where the infimum is over all estimators $\widehat{\eta}$ of the change point location and $\eta(P)$ denotes the change point location of $P \in \mathcal{Q}^T$.	
\end{lemma}

\Cref{lemma-error-opt} shows the minimax lower bound on the localisation error is of order $\sigma^2\kappa^{-2}$.  The localisation errors achieved by the estimators from Algorithms~\ref{algorithm:PDP} and \ref{algorithm:WBS} can both be nearly optimal, off by logarithmic factors, under suitable conditions.

\subsection{Conclusions}

The univariate piecewise constant change point detection and localisation are the blueprints for more complicated situations.  The optimal localisation rate is achievable under the nearly minimax optimal signal-to-noise ratio regime.  Recalling \Cref{sec-what-we-will-cover}, in the univariate mean change point problem, all three optimality goals are achieved, saving for logarithmic factors.  However, this phenomenon is not always true and we will see later.

We remark on some comparisons between Algorithms~\ref{algorithm:PDP} and \ref{algorithm:WBS}.
\begin{itemize}
	\item There is only one tuning parameter $\lambda$ in \Cref{algorithm:PDP}, but there are in fact two in \Cref{algorithm:WBS}, $\tau$ and $C_R$ (involved in the upper bound on the random interval lengths, see \Cref{thm-wbs-uni}).  Since all tuning parameters need to be specified in practice, \Cref{algorithm:PDP} is superior than \Cref{algorithm:WBS} in this aspect.
	\item The worst-case computational costs for Algorithms~\ref{algorithm:PDP} and \ref{algorithm:WBS} are of order $O(T^2)$ and $O(T^3)$, respectively.
	\item Beyond univariate mean change point problems, one key in the construction of the optimisation problem in \Cref{algorithm:PDP} is a proper choice of the cost function $H(\cdot)$ in \eqref{eq-G}.  This is not a problem in many cases, but might be a problem in nonparametric cases, which we will discuss in \Cref{sec-np-conclusions}.  In this case, \Cref{algorithm:WBS} may enjoy some flexibility in constructing a corresponding CUSUM statistic according to the choice of $\|\cdot\|_*$ in the model assumption.
\end{itemize}


In this section so far, we have only studied the situations where the data are independent and identically distributed between two consecutive change points.  There have been some works on the possible relaxations.
	\begin{itemize}
	\item Temporal dependence.  There are two popular ways to impose temporal dependence: one is through the noise sequence and the other is to assume the data are from some time series models.  

	As for the former, if one aims for fixed sample results, the results can easily be extended by using concentration inequalities developed for dependent data \citep[e.g.][]{delyon2009exponential}.  This method is adopted in \cite{padilla2019change}, among others.  Alternatively, if only asymptotic results are required, then one can just assume that the noise sequence is uncorrelated.  A scaled sum of the noise sequence can be shown to follow a Brownian bridge in the limit \citep[e.g.][]{aue2008testing, Lavielle1999}.
	
	As for the latter, the existing literature handles correlated data includes \cite{wang2020detecting}, \cite{wang2019localizing},\cite{dette2018multiscale}, \cite{akashi2018change}, \cite{dette2018change} and \cite{aue2009break}, among others.

	\item Robust estimators.  If the data between two consecutive change points are not necessarily identically distributed, then certain forms of robust estimation is required.  Works along this line include \cite{fearnhead2019changepoint}, \cite{pein2015heterogeneous} and \cite{yu2019robust}, among others.
	\end{itemize}

\section{Extension 1: Piecewise polynomials}\label{sec-piecewise-polynomial}

\subsection{Overview}

The studies on the univariate piecewise constant change point detection problem lays the foundation for studying more complicated problems.  One direction to generalise the results we discussed in \Cref{sec-detec-local} is to piecewise polynomials with any arbitrary but fixed orders.  In this section, we are concerned with the model, for each $t \in \{1, \ldots, T\}$,
	\begin{equation}\label{eq-y-intro}
		X_t = \theta_t + \varepsilon_t = f(t/T) + \varepsilon_t, 
	\end{equation}
	where $f(\cdot)$ is an unknown function belonging to the class $\mathcal{F}^{r, K}$, defined as
	\begin{align}
		\mathcal{F}^{r, K} = \Big\{f(\cdot): [0, 1] \to \mathbb{R}: & f \mbox{ has } K+1 \mbox{ pieces and each piece is a right-continuous }\nonumber \\
		& \mbox{with left limit polynomial of order at most } r\Big\}, \label{eq-def-f-r-n-k}
	\end{align}
		
The first task we have is to quantify the changes at every change point.	

\begin{definition}\label{def-jump-size}
Let $f(\cdot) \in \mathcal{F}^{r, K}$, $\{s_k\}_{k = 1}^K$ be the collection of all the change points of $f(\cdot)$, and $s_0 = 0$, $s_{K+1} = 1$.  For any $k \in \{1, \ldots, K\}$, let $f_{[s_{k-1}, s_{k+1})}(\cdot): [s_{k-1}, s_{k+1}) \to \mathbb{R}$ be the restriction of $f(\cdot): [0, 1] \to \mathbb{R}$ on $[s_{k-1}, s_{k+1})$.  Define the reparameterisation of $f_{[s_{k-1}, s_{k+1})}(\cdot)$ as
	\begin{equation}\label{eq-repara}
		f(x) = \begin{cases}
 			\sum_{l = 0}^r a_l (x-s_k)^l, & x \in [s_{k-1}, s_k), \\
 			\sum_{l = 0}^r b_l (x-s_k)^l, & x \in [s_k, s_{k+1}),
 		\end{cases}
	\end{equation}
	where $\{a_l, b_l\}_{l=0}^r \subset \mathbb{R}$.  Define the jump associated with the change poin $s_k$ as 
	\[
		\kappa_k = |a_{r_k} - b_{r_k}| > 0,
	\]
	where 
	\begin{equation}\label{eq-rk-defi}
		r_k = \min\{l = 0, \ldots, r: \, a_l \neq b_l\}.
	\end{equation}
\end{definition}

\Cref{def-jump-size} provides the definition of the jump size we concern in this problem.  In the following, we lay out the counterpart of \Cref{assume:change-uni-mean} in the piecewise polynomial case.

\begin{assumption}\label{assume-model-poly}
Assume that the data $\{X_t\}_{t = 1}^T$ are generated from \eqref{eq-y-intro}, where $f(\cdot)$ belongs to $\mathcal{F}^{r, K}$ defined in \eqref{eq-def-f-r-n-k} and $\varepsilon_i$'s are independent zero mean sub-Gaussian random variables with $\max_{i = 1}^n \|\varepsilon_i\|_{\psi_2} \leq \sigma^2$. 

Let $\theta = (\theta_t)_{t = 1}^T$, with $\theta_t = f(t/T)$, be the discretised $f(\cdot)$ on the grid of $\{1/T, 2/T, \ldots, 1\}$.  We denote the collection of all change points of $\theta$ to be $\{\eta_1, \ldots, \eta_K\}$, satisfying
	\[
		\Delta = \min_{k \in \{1, \ldots, K+1\}} (\eta_k - \eta_{k-1}) > 0
	\]
	where $\eta_0 = 1$ and $\eta_{K+1} = T+1$.
	
In addition, for any $k \in \{1, \ldots, K\}$, let 
	\[
		\kappa = \min_{k = 1, \ldots, K} \kappa_k > 0,
	\]
	where $\kappa_k$ is defined in \Cref{def-jump-size}.
\end{assumption}

Comparing to \Cref{assume:change-uni-mean}, in \Cref{assume-model-poly} we can see that the underlying signals are allowed to be any arbitrary but fixed order of polynomials, instead of just constants.  Besides this apparent difference, we would like to highlight a few more.
	\begin{itemize}
		\item With the sample size $T$, the time scales in Assumptions~\ref{assume:change-uni-mean} and \ref{assume-model-poly} are $O(1)$ and $O(1/T)$, respectively.  These two are in fact equivalent, but we adopt the two different scales to follow the suit in the existing literature.  Estimating piecewise polynomial signals has a rich body of literature, including \cite{shen2020phase}, \cite{mammen1997locally}, \cite{tibshirani2014adaptive},\cite{rudin1992nonlinear}, \cite{zhang2002risk} and \cite{chatterjee2015risk}, among others.  
		\item The jumps in \Cref{assume:change-uni-mean} are characterised by the mean changes, which are natural due to the piecewise-constant features.  For two different at-most-order-$r$ polynomials, they are specified by two coefficient vectors.  There are different ways to measure the difference between two different coefficient vectors.  We characterise the distance in \Cref{assume-model-poly}(c) -- this provides the sharpest localisation rates.
	\end{itemize}

\subsection{Consistent localisation}

Recall that the penalised estimator we studied in \Cref{sec-uni-mean-ell-0} is a penalised sum of residuals.  The residuals are defined to be the residuals after projecting data onto the space spanned by the all one vector.  To be specific, for the interval $I$ and its corresponding data vector $X_I = (X_i, i \in I)^{\top}$, the projection matrix is defined to be $\Pi_{I, 0} = \Pi_{I} = \mathbbm{1}_{|I|} (\mathbbm{1}_{|I|}^{\top}\mathbbm{1}_{|I|})^{-1}\mathbbm{1}_{|I|}$.  We add the extra subscript 0 in $\Pi_I$, since constants are order-0 polynomials.  When we move from piecewise constant signals to piecewise polynomial signals, one can generalise the projection matrix correspondingly.  
	
Let $I = [s, e] \subset \{1, \ldots, T\}$, $r \in \mathbb{N}$ and
	\begin{equation}\label{eq-U-I-D-definition}
		U_{I, r} = \left(\begin{array}{cccc}
				1 & s/T & \cdots & (s/T)^r \\
				\vdots & \vdots & \vdots & \vdots \\
				1 & e/T & \cdots & (e/T)^r  
			\end{array}\right) \in \mathbb{R}^{(e-s+1) \times (r+1)}.
	\end{equation}
	We define 
	\begin{equation}\label{eq-Pi-I-D-definition}
		\Pi_{I, r} = U_{I, r}(U_{I, r}^{\top}U_{I, r})^{-1}U_{I, r}^{\top} 
	\end{equation}	
	to be the order-$r$ polynomial projection matrix.  The change point estimators are the output of \Cref{algorithm:PDP} with
	\begin{equation}\label{eq-H-poly}
		H(I) = \|X_I - \Pi_{I, r}X_I\|^2.
	\end{equation}
	The theoretical guarantees of the outputs are given below.

\begin{theorem}\label{thm-main-poly}
Let $\{X_t\}_{t=1}^T$ satisfy \Cref{assume-model-poly}.  In addition, assume there exists a large enough constant $C_{\mathrm{SNR}} > 0$ and any $\xi > 0$, such that
	\begin{equation}\label{eq-snr-assum-poly}
		\min_{k = 1, \ldots, K}\frac{\kappa_k^2 \Delta^{2 r_k + 1}}{\sigma^2 T^{2r_k}} \geq C_{\mathrm{SNR}} K \log^{1+\xi}(T).
	\end{equation}
	With $\lambda = C_{\lambda} K\sigma^2\log(T)$, let $\{\widehat{\eta}_k\}_{k = 1, \ldots, \widehat{K}}$ be the collection of change point estimators from \Cref{algorithm:PDP}, with $H(\cdot)$ defined in \eqref{eq-H-poly} and $C_{\lambda} > 0$ being an absolute constant, satisfy that 
	\begin{equation}\label{eq-thm-1-poly}
		\mathbb{P}\left\{\widehat{K} = K, \, \forall k \in \{1, \ldots, K\},\, \left|\widehat{\eta}_k - \eta_k\right| \leq \left[\frac{C_{\epsilon} K T^{2r_k}\sigma^2 \log(T)}{\kappa_k^2}\right]^{1/(2r_k+1)} \right\} > 1 - T^{-c},
	\end{equation}
	where $c, C_{\epsilon} > 0$ are absolute constant.
\end{theorem}		

The assumption \eqref{eq-snr-assum-poly} and the localisation error in \eqref{eq-thm-1-poly} show that the penalised estimator is consistent.  To be specific,
	\begin{align*}
		& \frac{1}{\Delta}\max_{k = 1, \ldots, K}\left[\frac{C_{\epsilon} K T^{2r_k}\sigma^2 \log(T)}{ \kappa_k^2}\right]^{1/(2r_k+1)}  \\
		& \hspace{2cm} \lesssim \max_{k = 1, \ldots, K} \left[\frac{K \sigma^2 \log(T) T^{2_{r_k}}}{\kappa_k^2 \Delta^{2r_k +1}}\right]^{1/(2r_k + 1)} \lesssim \log^{-1/(2r_k+1)}(T) \to 0,
	\end{align*}
	as $T$ grows unbounded.  
	
Note that in \eqref{eq-thm-1-poly}, the larger $r_k$ is, the larger the localisation error is.  This explains our choice of distribution difference in \Cref{def-jump-size}.  Choosing the smallest order with different coefficients yields sharpest localisation errors.

\subsection{Optimal localisation}

The signal-to-noise ratio and the localisation error lower bounds are presented in Lemmas~\ref{lem-lb-1-poly} and \ref{lem-lb-2-poly}, respectively.  

\begin{lemma}\label{lem-lb-1-poly}
Under \Cref{assume-model-poly}, assume that there exists one and only one change point and $d_1 = r$.  Let $P_{\kappa, \Delta, \sigma, r, T}$ denote the joint distribution of the data.  For a small enough $c_1 > 0$, consider the class
	\[
		\mathcal{P}^T = \left\{P_{\kappa, \Delta, \sigma, r, T}: \, \Delta = \min \left\{\Bigg\lfloor \left(\frac{c_1 T^{2r}}{\kappa^2 \sigma^{-2} }\right)^{1/(2r+1)} \Bigg\rfloor,\, T/3 \right\} \right\}.
	\]	
	Then we have
	\[
		\inf_{\hat{\eta}} \sup_{P \in \mathcal{P}^T} \mathbb{E}_P(|\hat{\eta} - \eta(P)|) \geq c,
	\]
	where $\eta(P)$ is the location of the change point for distribution $P$, the minimum is taken over all the measurable functions of the data and $0 < c < 1$ is an absolute constant depending on $c_1$.
\end{lemma}

\begin{lemma}\label{lem-lb-2-poly}
Under \Cref{assume-model-poly}, assume that there exists one and only one change point and $d_1 = r$.  Let $P_{\kappa, \Delta, \sigma, r, T}$ denote the joint distribution of the data.  Consider the class
	\[
		\mathcal{Q}^T = \left\{P_{\kappa, \Delta, \sigma, r, T}: \, \Delta < T/2,  \, \kappa^2 \Delta^{2r+1} \geq \sigma^2 T^{2r} \zeta_T\right\},
	\]
	for any diverging sequence $\{\zeta_T\}$.  Then for all $T$ large enough, it holds
	\[
		\inf_{\hat{\eta}} \sup_{P \in \mathcal{Q}^T} \mathbb{E}_P(|\hat{\eta} - \eta(P)|) \geq \max\left\{1, \, \left[\frac{c \sigma^2}{T \kappa^2}\right]^{1/(2r+1)}\right\},
	\]
	where $\eta(P)$ is the location of the change point for distribution $P$, the minimum is taken over all the measurable functions of the data and $0 < c < 1$ is an absolute constant.
\end{lemma}

We can see from Lemmas~\ref{lem-lb-1-poly} and \ref{lem-lb-2-poly} that both the signal-to-noise ratio condition and the localisation error we have in \Cref{thm-main-poly} are off by a logarithmic factor and another factor of $K$.  Since $K$ is allowed to diverge, \Cref{thm-main-poly} is sub-optimal in detection and localisation.  We conjecture that this sub-optimality is due to an artefact of the proof and the optimisation problem we study in \Cref{thm-main-poly} should have been nearly-optimal only off by a logarithmic factor.  Having said this, in order to improve, we present a refinement step.

\begin{theorem}\label{thm-refine}
Under all the assumptions in \Cref{thm-main-poly}, let $\{\nu_k\}_{k = 1}^K$ satisfy
	\begin{equation}\label{eq-refine-thm-loc-cond}
		T \max_{k = 1, \ldots, K}|\nu_k - \eta_k| \leq \Delta/5.
	\end{equation}
	For each $k \in \{1, \ldots, K\}$, define
	\[
		s_k = T\nu_{k-1}/2 + T\nu_k/2, \quad e_k = T\nu_k/2 + T\nu_{k+1}/2 \quad \mbox{and}\quad I_k = (s_k, e_k),
	\]
	with $\nu_0 = 1/T$ and $\nu_{K+1} = 1+1/T$.  For $k \in \{1, \ldots, K\}$, we let
	\[
		\widetilde{\eta}_k = \min_{t \in I_k} \left(\|X_{[s_k, t)} - \Pi_{[s_k, t), r}X_{[s_k, t)}\|^2 + \|X_{[t, e_k)} - \Pi_{[t, e_k), r}X_{[t, e_k)}\|^2\right),
	\]
	where $\Pi_{\cdot, r}$ is defined in \eqref{eq-U-I-D-definition} and \eqref{eq-Pi-I-D-definition}.  Then we have
	\[
		\mathbb{P}\left\{\forall k \in \{1, \ldots, K\},\, \left|\widetilde{\eta}_k - \eta_k\right| \leq \left[\frac{C_{\epsilon}\sigma^2 T^{2r_k}\log(T)}{\kappa_k^2}\right]^{1/(2r_k+1)} \right\} > 1 - T^{-c},
	\]
	where $c, C_{\epsilon} > 0$ are absolute constant.
\end{theorem}	

\Cref{thm-refine} guarantees that, under \eqref{eq-refine-thm-loc-cond}, one can achieve nearly optimal localisation errors.  Since \Cref{thm-main-poly} shows the outputs of \Cref{algorithm:PDP} with \eqref{eq-H-poly} satisfy \eqref{eq-refine-thm-loc-cond} with large probability, \Cref{thm-refine} can be used as a second step to refine the outputs of \Cref{algorithm:PDP}.

As for the signal-to-noise ratio, there is still a gap of order $K$ between the detection upper and lower bounds, therefore we are not able to fulfil the three goals listed in \Cref{sec-what-we-will-cover}.  We conjecture that this gap is due to a loose upper bound and the loose upper bound is due to a loose control of some cross terms.  We explain it below.

\begin{lemma}\label{lem-diff}
Let $I_1$ and $I_2$ denote any two disjoint intervals of $\{1, \ldots, T\}$ and $I = I_1 \cup I_2$.  For any sequences $\{X_i\}_{i = 1/T, 2/T, \ldots, 1} \subset \mathbb{R}$, it holds that 
	\begin{align*}
		\|X_I - \Pi_{I, r} X_I\|^2  = \|X_{I_1} - \Pi_{I_1, r}X_{I_1}\|^2 + \|X_{I_2} - \Pi_{I_2, r} X_{I_2}\|^2 + Q(I_1, I_2, \{X\}, r),
	\end{align*}
	where 
	\begin{align*}
		& Q(I_1, I_2, \{X\}, r) \\
		= & \{X_{I_1}^{\top}U_{I_1, r}(U_{I_1, r}^{\top}U_{I_1, r})^{-1} - X_{I_2}^{\top}U_{I_2, r}(U_{I_2, r}^{\top}U_{I_2, r})^{-1}\} \{(U_{I_1, r}^{\top}U_{I_1, r})^{-1} + (U_{I_2, r}^{\top}U_{I_2, r})^{-1}\}^{-1} \\
		& \hspace{2cm} \times \{(U_{I_1, r}^{\top}U_{I_1, r})^{-1} U_{I_1, r}^{\top}X_{I_1} - (U_{I_2, r}^{\top}U_{I_2, r})^{-1} U_{I_2, r}^{\top}X_{I_2}\}
	\end{align*}
	and $U_{\cdot, r}$ is defined in \eqref{eq-U-I-D-definition}.
\end{lemma}

The key to provide localisation errors is to lower bound the cross term $Q(I_1, I_2, \{X\}, r)$, when partitioning $I$ into $I_1$ and $I_2$ provides a good estimator of a change point; and to upper bound $Q(I_1, I_2, \{X\}, r)$, if partitioning $I$ into $I_1$ and $I_2$ leads to over-partitioning.  

Note that when $r = 0$, i.e.~in the piecewise constant case,
	\[
		Q(I_1, I_2, \{\mathbb{E}(X)\}, 0) = \frac{|I_1||I_2|}{|I_1| + |I_2|} \left(|I_1|^{-1}\sum_{i \in I_1} \mathbb{E}(X_i) - |I_2|^{-1}\sum_{i \in I_2} \mathbb{E}(X_i)\right)^2.
	\]
	In addition, it holds that
	\[
		\frac{\min\{|I_1|, \, |I_2|\}}{2} = \frac{|I_1||I_2|}{2 \max\{|I_1|, \, |I_2|\}} \leq \frac{|I_1||I_2|}{|I_1| + |I_2|} \leq \min\{|I_1|, \, |I_2|\}.
	\]
	If $I_1 = [s, \eta)$ and $I_2 = [\eta, e)$, and if $\eta$ is the only true change point in $I_1 \cup I_2$, then we have
	\[
		Q(I_1, I_2, \{\mathbb{E}(X)\}, 0) \asymp \min\{|I_1|, \, |I_2|\} \kappa^2.
	\]	
	
As for general $r$ and for the case discussed above, we have that
	\begin{equation}\label{eq-Q-lack-upper}
		Q(I_1, I_2, \{\mathbb{E}(X)\}, r) \geq \frac{\min\{|I_1|, \, |I_2|\}}{2} \kappa^2,
	\end{equation}
	but lacks an upper bound of the same order.  The lack of such an upper bound directly resulted in the term $K$ in \eqref{eq-snr-assum-poly}.  We remark that, when $r = 1$, a similar estimator is studied in \cite{fearnhead2019detecting}, where $K$ is assumed to be an absolute constant.  In general, an ideal solution to match the three goals in \Cref{sec-what-we-will-cover} is yet known.  We conjecture the ideal solution is reachable if one can provide an upper bound of order $\min\{|I_1|, \, |I_2|\} \kappa^2$ in \eqref{eq-Q-lack-upper}.	

\subsection{Conclusions}\label{sec-poly-sec-conclusions}

The parting words of this section are about the refinement step we discussed in \Cref{thm-refine}.  The refinement idea will appear again in more complicated situations, for instance the high-dimensional graphon in \Cref{theorem:localization 1}, the high-dimensional linear regression in \Cref{cor-lr-high-dim} and \cite{wang2019statistically}, in the high-dimensional vector autoregressive models in \cite{wang2019localizing}.  

The motivation of adopting an additional step is that directly localising multiple change points may lead to a consistent but not necessarily optimal localisation error rates.  An additional step is useful if at least one of the following situations holds.
	\begin{itemize}
	\item The sub-optimality of the directly localising multiple change points is due to the multiple change points.  This is the case we have in \Cref{thm-main-poly}.  The additional step works in the interval contains, with large probability, one and only one true change point, and therefore improves the rate with respect to the number of change points.
	\item The sub-optimality of the directly localising multiple change points is due to the choice of estimators of the underlying distributions, and there exists better estimators of the underlying distributions.  This is the case we have in Theorems~\ref{theorem:localization 1} and \ref{cor-lr-high-dim}.  

		Although we emphasised at the beginning that this survey is not covering in-depth results and discussions on the computational aspect, this is a serious issue especially in the more complicated data type scenarios.  In order to save computational costs, in the multiple change points case, one may first use some computationally-cheaper estimation in either the $H(\cdot)$ function in \eqref{eq-G} or the CUSUM statistic defined in \Cref{def-cusum}.  Then in the refinement step, since the optimisation in each working interval is independent and there is only one change point in each interval, one may want to adopt some computationally-more-expensive estimator to yield better estimation.
	\end{itemize}

\section{Extension 2: High-dimensional problems}

As we have mentioned in \Cref{sec-intro}, the change point analysis is by no means restricted to detecting changes in a sequence of univariate data.  In this section, we consider three high-dimensional extensions, the data of which are sequences of high-dimensional vectors, high-dimensional matrices and high-dimensional regression coefficients, respectively.   For each of these three scenarios, we use one type of change point problem to illustrate, and we conclude this section with other problems studied in the existing literature, in addition to some insights on how the high-dimensionality affects the difficulties of the problems.

\subsection{Covariance changes}\label{sec-cov-change}

The first high-dimensional case we consider is a sequence of high-dimensional random vectors, the covariances of which are piecewise constant.  The model is detailed below.

\begin{assumption}\label{assume:model-cov}
Let $\{X_t\}_{t=1}^T \subset \mathbb{R}^p$ be independent, zero mean random vectors such that $\mathbb{E}(X_tX_t^{\top}) = \Sigma_t$ and $\|X_t\|_{\psi_2}\le \sigma$ for all $t = 1, \ldots, T$, where $\sigma >0$. Let  $\{\eta_0,\ldots,\eta_{K+1} \} \subset \{1,\ldots, T+1\}$ be a strictly increasing subsequence of change points such that  $\eta_0=1$,  $\eta_{K+1}=T+1$ and  $\Sigma_t \neq \Sigma_{t-1}$ if and only if  $t \in \{\eta_1,\ldots,\eta_{K}\}$.  The minimal spacing between jumps is defined to be
		\[
			\Delta = \min_{k = 1, \ldots, K+1} \{\eta_k-\eta_{k-1}\} > 0,
		\]
		and the magnitude of changes is 
		\[
			\kappa = \min_{k = 1, \ldots, K} \kappa_k = \min_{k = 1, \ldots, K}\|\Sigma_{\eta_k} -\Sigma_{\eta_{k}-1} \|_{\mathrm{op}} > 0.
		\]
\end{assumption}

Note that the $\|\cdot\|_{\psi_2}$ is the Orlicz-$\psi_2$-norm or the sub-Gaussian norm, defined as, for any random vector $X \in \mathbb{R}^p$,
	\[
		\|X\|_{\psi_2} = \sup_{v \in \mathbb{R}^p, \, \|v\| = 1} \|X^{\top}v\|_{\psi_2}.
	\]

The parameters $\sigma$ and $\kappa$ reflect the magnitudes of noise and signal, respectively.  However they are not variation independent, as they
    satisfy the inequality $\kappa \leq  \sigma^2/4$, due to the following derivations:
    \begin{align*}
		\kappa \leq \max_{k=1}^{K} \|\Sigma_{\eta_k} - \Sigma_{\eta_{k}-1}\|_{\mathrm{op}} \leq 2 \max_{t=1}^T\|\Sigma_t\|_{\mathrm{op}} = 2\max_{t=1}^T \sup_{v \in \mathbb{R}^p, \, \|v\| = 1} \mathbb{E}\bigl[(v^{\top}X_t)^2\bigr] \leq 4\max_{t=1}^T\|X_t\|_{\psi_2}^2 \leq 4\sigma^2.
	\end{align*}
	This is a trademark of the covariance change point problems and provides extra difficulties as opposed to mean change point problems -- the larger the jumps are, the larger the variances are.
	
To thoroughly understand the difficulties of high-dimensional covariance change point problems, we provide the minimax lower bounds on the detection and localisation errors, which are collected in Lemmas~\ref{lemma:lower bound 1-cov} and \ref{lemma:lower bound 2-cov}, respectively.	
	
\begin{lemma}[Lemma 3 in \citealp{wang2017optimal}] \label{lemma:lower bound 1-cov}
Under \Cref{assume:model-cov}, assume that there is one and only one change point.  Let $P_{\kappa, \Delta, \sigma, p, T}$ denote the joint distribution of the data.  Consider the class of distributions 
	\[
		\mathcal P^T = \left\{P_{\kappa,\Delta, \sigma, p, T} \colon \Delta \leq \min \left\{ \frac{2\sigma^4p}{33\kappa^2}, \, T/3 \right\},   \, \kappa\le \sigma^2/4 \right\}.
	\]
	We have that,
	\[
    	\inf_{\widehat{\eta}} \sup_{P\in \mathcal P^T} \mathbb{E}_P(|\widehat{\eta} -\eta(P)|) \geq \Delta/2,
	\]
	where $\eta(P)$ is the location of the change point of distribution $P$ and the infimum is over all estimators of the change point. 
\end{lemma}	

\begin{lemma}[Lemma 4 in \citealp{wang2017optimal}] \label{lemma:lower bound 2-cov}
Under \Cref{assume:model-cov}, assume that there is one and only one change point.  Let $P_{\kappa, \Delta, \sigma, p, T}$ denote the joint distribution of the data.  Consider the class of distributions 
	\[
		\mathcal Q^T = \left\{ \mathcal{P}_{\kappa,\Delta,\sigma, p, T} \colon \,  \Delta \kappa^2 \ge  p  \log(T)\sigma^4, \, \kappa\le \sigma^2/4, 4 \leq \Delta \leq 4/5 (T-1) \right\}.
	\]
	Then,
	\[
	\inf_{\widehat{\eta}} \sup_{P\in \mathcal Q^T} \mathbb{E}_P(|\widehat{\eta} -\eta(P)|) \geq \frac{\sigma^4}{20 \kappa^2},
	\]
where $\eta(P)$ is the location of the change point of distribution $P$ and the infimum is over all estimators of the change point. 
\end{lemma}

Lemmas~\ref{lemma:lower bound 1-cov} and \ref{lemma:lower bound 2-cov} show that in the low signal-to-noise ratio regime $\kappa^2 \Delta \lesssim p\sigma^4$, no algorithm is guaranteed to produce consistent change point estimators.  The optimal localisation error in the high signal-to-noise ratio regime $\kappa^2 \Delta \gtrsim p \sigma^4 \log(T)$ is $\sigma^4 \kappa^{-2}$. \\
	
A scan-statistics-based algorithm is able to achieve the near-optimality in the sense of detection boundary and localisation error. 
The detailed algorithm is given in \Cref{algorithm:WBSRP}, with a subroutine specified in \Cref{algorithm:PC} and theoretical guarantees available in \Cref{thm:wbsrp}.  The quantity $\widetilde{Y}^{s_m, e_m}_t(u_m)$ is a CUSUM statistic defined in \Cref{def-cusum}.
	
\begin{algorithm}[!ht]
\begin{algorithmic}
	\INPUT $\{X_t\}_{t=1}^T$, $\{ (\alpha_m ,\beta_m) \}_{m=1}^M$
	 \For{$m = 1, \ldots, M$}
	 \If {$\beta_m-\alpha_m> 2p\log(T) +1 $}
		\State $d_m \leftarrow \argmax_{\lceil \alpha_m +  p\log(T) \rceil \leq t \leq \lfloor \beta_m -  p\log(T) \rfloor} \|\widetilde{S}^{\alpha_m, \beta_m}_t \|_{\mathrm{op}}$
		\State $u_m \leftarrow \argmax_{\|v\| =1} \bigl| v^{\top}\widetilde{S}^{\alpha_m, \beta_m}_{d_m}v|$
		\Else
		\State $u_m\leftarrow0$
	\EndIf
	\EndFor	
	\OUTPUT $\{ u_m\}_{m=1}^M$.
\caption{Principal Component Estimation $\mathrm{PC}(\{X_t\}_{t=1}^T, \{(\alpha_m ,\beta_m)\}_{m=1}^M)$}
\label{algorithm:PC}
\end{algorithmic}
\end{algorithm}

\begin{algorithm}[!ht]
\begin{algorithmic}
	\INPUT Two independent samples, $\{W_t\}_{t=1}^{T}$ and $\{X_t\}_{t=1}^{T}$, and the threshold parameter $\tau>0$.
	\State  $\{ u_m\}_{m=1}^M \leftarrow PC(\{W_t\}_{t=1}^T, \{(\alpha_m ,\beta_m)\}_{m=1}^M  )$
	\For{$t \in  \{s, \ldots, e\}$}
		\For{$m = 1, \ldots, M$}
			\State $Y_t(u_m) \leftarrow \bigl(u_m^{\top} X_t\bigr)^2$
		\EndFor	
	\EndFor		
	\For{$m = 1, \ldots, M$}  
		\State $(s_m , e_m ) \leftarrow [s,e]\cap [\alpha_m,\beta_m]$
		\If{$e_m - s_m \geq 2\log(T) + 1$}
			\State $b_{m} \leftarrow \argmax_{s_m +\log(T) \leq t \leq e_m -\log(T)}  | \widetilde Y^{s_m,e_m}_{t} (u_{m})|$  
			\State $a_m \leftarrow \bigl| \widetilde Y^{s_m,e_m}_{b_{m}}  ( u_{m})\bigr|$ 
		\Else 
			\State $a_m \leftarrow -1$	
		\EndIf
	\EndFor	
	\State $m^* \leftarrow \argmax_{m = 1, \ldots, M} a_{m}$
	\If{$a_{m^*} > \tau$}
		\State add $b_{m^*}$ to the set of estimated change points
		\State WBSIP$(\{X_t, W_t\}_{t = 1}^T, (s, b_{m*}),\{ (\alpha_m,\beta_m)\}_{m=1}^M, \tau )$
		\State WBSIP $(\{X_t, W_t\}_{t = 1}^T, (b_{m*}+1,e),\{ (\alpha_m,\beta_m)\}_{m=1}^M,\tau ) $
	\EndIf
	\OUTPUT The set of estimated change points.
\caption{Wild Binary Segmentation through Independent Projection. WBSIP$(\{X_t, W_t\}_{t = 1}^T, (s, e),$ $\{ (\alpha_m,\beta_m)\}_{m=1}^M, \tau $)}
\label{algorithm:WBSRP}
\end{algorithmic}
\end{algorithm} 
	
The requirement of two independent samples in \Cref{algorithm:WBSRP} can be achieved by splitting data into even and odd indices subsets.	
	
\begin{theorem}[Theorem 2 in \citealp{wang2017optimal}] \label{thm:wbsrp}
Let \Cref{assume:model-cov} hold and let $\{(\alpha_m,\  \beta_m) \}_{m=1}^M\subset (0, T)$ be a collection of intervals whose endpoints are drawn independently and uniformly from $\{1,\ldots, T\}$ and such that $\max_{1\le m \le M} (\beta_m -\alpha_m)\le C_R \Delta$  for an absolute constant $C_R>0$.  In addition, assume that for any $\xi > 0$, there exists a sufficiently large absolute constant $C >0$ such that 
	\[
		\Delta \kappa^2\ge C p\log^{1+\xi}(T) \sigma^{4}.
	\]
	
Suppose there exist sufficiently small constant  $c_2>0$ and sufficiently large constant  $c_3 > 0$   such that the input parameter   $\tau$  satisfy
	\[
		c_3 \sigma^2 \sqrt {\log(T)}< \tau <c_2\kappa \sqrt {\Delta}.	
	\]
	Then the collection of the estimated change points $\{\widehat \eta_k\}_{k=1}^{\widehat K}$	returned by \Cref{algorithm:WBSRP} with input parameters of  $(0, T)$, $\{(\alpha_m,\beta_m)\}_{m=1}^M$ and $\tau$, satisfies	
		\begin{align}
			 & \mathbb{P}\Bigl\{\widehat K =K \quad \mbox{and} \quad   |\eta_k-\widehat \eta_k| \leq C_1 \sigma^4\log(T) \kappa_k^{-2}, \, \forall  k \in \{1, \ldots, K\} \Bigr\} \nonumber \\
			& \hspace{2cm}\geq  1-4 M  T^{2-c}  - 4\times 9^p T^ {3-cp}  -\exp\bigl\{\log(T/\Delta)- M\Delta/(4C_R T)\bigr\}, \nonumber
	 	\end{align}
	for some absolute constants $c > 3$ and $C_1 > 0$.
\end{theorem}

Provided that $\log(T/\Delta) \lesssim M\Delta/T$, we see from \Cref{thm:wbsrp} that \Cref{algorithm:WBSRP} is nearly optimal in terms of both detection and localisation, despite that the dimension $p$ is allowed to grow unbounded as the sample size diverges.  This means all three goals we listed in \Cref{sec-what-we-will-cover} are achieved.  It might come as a surprise that the optimal localisation error rate is not a function of the dimension $p$, and we are actually able to achieve it under the minimal conditions.  We will come back to discuss this phenomenon in \Cref{sec-high-dim-conclusion}, together with more high-dimensional cases.  

Finally, the covariance change point analysis has also been studied in different settings over the years, including \cite{InclanTiao1994}, \cite{GombayEtal1996}, \cite{DetteEtal2018}, \cite{avanesov2016change}, \cite{BirkeDette2005} and \cite{aue2009break}, among others.

\subsection{Graphon changes}\label{sec-graphon}

\subsubsection{Overview}
Instead of obtaining a random vector at every time, the random objects obtained can be in the form of random matrices.  With the surging of network data, we use a dynamic networks model as an example and study the graphon changes in this subsection.  The detailed model assumptions are collected in \Cref{assume:model-network}, with a general definition in \Cref{def-inhomo-ber-net}.

\begin{definition}[Inhomogeneous Bernoulli networks]\label{def-inhomo-ber-net}
A network with node set $\{1, \ldots, n\}$  is an inhomogeneous Bernoulli network if its adjacency matrix $A \in \mathbb{R}^{n\times n}$ satisfies
		\[
		A_{ij} = A_{ji} = \begin{cases}
			1, & \mbox{nodes $i$ and $j$ are connected by an edge},\\
			0, & \mbox{otherwise};
		\end{cases}
		\]
		and $\{A_{ij}, i < j\}$ are independent Bernoulli random variables with $\mathbb{E}(A_{ij}) = \Theta_{ij}$. 
\end{definition}

\begin{assumption}\label{assume:model-network}
	Let $\{X_t\}_{t = 1}^T \subset \mathbb{R}^{n\times n}$ be a collection of adjacency matrices of independent inhomogeneous Bernoulli networks with means $\{\Theta_t\}_{t = 1}^T$ satisfying the following properties.

The sparsity parameter  
		\[
			\rho = \max_{t = 1, \ldots, T} \|\Theta_t\|_{\infty}
		\]
	is such that $\rho n  \ge \log(n)$, where $\|\cdot\|_{\infty}$ denotes the entrywise maximum norm of a matrix.

There exists a sequence $\eta_0 < \eta_1 < \ldots < \eta_{K+1}$ of time points, called change points, with $\eta_0=1$ and $\eta_{K+1}=T+1$, such that 
	\[	
		\Theta_t \neq \Theta_{t-1}, \quad \mbox{if and only if} \quad t \in \{\eta_k\}_{k = 1}^K.
	\]

The minimal spacing between two consecutive change points satisfies 
	\[
	\min_{k = 1, \ldots, K+1} \{\eta_k-\eta_{k-1}\} = \Delta > 0.
	\]

The magnitudes of the changes in the data generating distribution are such such that
	\[
		\|\Theta (\eta_{k} ) -\Theta (\eta_{k}-1)\|_{\mathrm{F}} = \kappa_k, \quad  k=1, \ldots, K,
	\]
	where $\|\cdot\|_{\mathrm{F}}$ denotes the Frobenius norm of a matrix.  Let
	\[
		\kappa_0 = \frac{\kappa}{n\rho } = \frac{\min_{k = 1, \ldots, K}\kappa_k}{n\rho }.
	\]
\end{assumption}

For dynamic networks, the data are a sequence of adjacency matrices and their distributions are determined by their graphons, i.e.~the expectations of the adjacency matrices, if we assume the networks are inhomogeneous Bernoulli networks defined in \Cref{def-inhomo-ber-net}.  Then in terms of characterising jumps, it would be natural to seek a certain matrix norm.  In \Cref{sec-cov-change}, the matrix operator norm is adopted, and in \Cref{assume:model-network}, the matrix Frobenius norm is summoned.  One could argue that as for network models, the Frobenius norm is able to capture a richer collection of changes.  In this survey, we would pay more attention on how different choices of norms affect the difficulty of the problems.  We will come back to this in \Cref{sec-high-dim-conclusion}.

The difficulty of the graphon change point analysis is explained in Lemmas~\ref{lemma:lower bound testing-network} and \ref{lem-3.3-lower}, on the minimax lower bounds on detection and localisation errors, respectively.

\begin{lemma}[Lemma 1 in \citealp{wang2018optimal}] \label{lemma:lower bound testing-network}
Let $\{X_t\}_{t=1}^T$ be a sequence of independent inhomogeneous Bernoulli networks satisfying \Cref{assume:model-network} with $K = 1$.  Let $P_{\kappa_0, \Delta, n, \rho}^T$ denote the corresponding joint distribution.  For $\zeta \leq 1/33$, consider the class of distributions 
	\[
		\mathcal{P}^T = \left\{ P^T_{\kappa_0,\Delta,n ,\rho}:  \Delta = \min \biggl\{\bigg\lfloor \frac{\log(T)\zeta }{n\rho\kappa_0^2} \bigg\rfloor, \, \lfloor T/3 \rfloor \biggr\},\, \rho \leq 1/2, \,\kappa_0 \leq 1\right\}.
	\]
	For each $P \in\mathcal{P}^T$, let $\eta(P) \in \{1,\ldots,T\}$ denote the location of the corresponding change point.  It holds that	
	\[
		\inf_{\widehat{\eta}} \sup_{P \in \mathcal{P}^T} \mathbb{E}_P(|\widehat{\eta} -\eta(P)|) \geq 3\Delta/4,
	\] 
	where the infimum is over all the possible estimators of the change point location.
\end{lemma}

\begin{lemma}[Lemma 2 in \citealp{wang2018optimal}]\label{lem-3.3-lower}
Let $\{X_t\}_{t=1}^T$ be a sequence of independent inhomogeneous Bernoulli networks satisfying \Cref{assume:model-network} with $K = 1$.  Let $P_{\kappa_0, \Delta, n, \rho}^T$ denote the corresponding joint distribution. Consider the class of distributions 
	\[
		\mathcal Q^T = \left \{P_{\kappa_0, \Delta, n, \rho}^T: \kappa_0 \le 1/2, \, \rho \leq 1/2   \right\} .
	\]
	It holds that
	\[
		\inf_{\widehat{\eta}} \sup_{P \in \mathcal Q^T} \mathbb{E}_P(|\widehat{\eta} - \eta(P)|)\geq \max \{c\kappa_0^{-2}n^{-2}\rho, \, 1/2\},
	\]
	where $\eta(P)$ denotes the location of the corresponding change point and the infimum is over all the possible estimators of the change point location.
\end{lemma}

Lemmas~\ref{lemma:lower bound testing-network} and \ref{lem-3.3-lower} show that in the low signal-to-noise ratio regime $\kappa_0 \sqrt{\Delta} \lesssim (n \rho)^{-1/2} \log^{1/2}(T)$, no algorithm is guaranteed to be consistent, and the localisation error is lower bounded by $\rho\kappa_0^{-2} n^{-2}$.  To match these lower bounds, we provide two sets of algorithms.  
	\begin{itemize}
	\item \Cref{algorithm:MWBS} is theoretically supported by \Cref{thm-1-network}, showing that there exists a computationally-efficient method providing consistent change point estimators, with a nearly optimal signal-to-noise ratio condition.
	\item \Cref{algorithm:RI}, with a subroutine in \Cref{algorithm:USVT}, is theoretically supported by \Cref{theorem:localization 1}, showing that under a stronger condition, there exists a method providing nearly optimal localisation errors. 
	\end{itemize}

\subsubsection{Consistent localisation}

We study a scan-statistics-based algorithm, using a network CUSUM statistic, i.e.~\Cref{def-cusum} with $X_t$'s being adjacency matrices.

\begin{algorithm}[!ht]
	\begin{algorithmic}
		\INPUT Two independent samples $\{X_t\}_{t=1}^{T}, \{W_t\}_{t=1}^{T} \in \mathbb{R}^{n\times n}$, $\tau_1$.
		\For{$m = 1, \ldots, M$}  
			\State $[s_m', e_m'] \leftarrow [s,e]\cap [\alpha_m,\beta_m]$
			\State $(s_m, e_m) \leftarrow [s_m' + 64^{-1} (e'_m - s'_m),e_m' - 64^{-1} (e_m' - s_m')] $  
			\If{$e_m - s_m \ge 1$}
				\State $b_{m} \leftarrow \arg\max_{t = s_m+1, \ldots, e_m-1}  (\widetilde{X}^{s_m, e_m}(t), \widetilde{W}^{s_m, e_m} (t))$
				\State $a_m \leftarrow (\widetilde{X}^{s_m, e_m}(t), \widetilde{W}^{s_m, e_m} (t))$
			\Else 
				\State $a_m \leftarrow -1$	
			\EndIf
		\EndFor
		\State $m^* \leftarrow \arg\max_{m = 1, \ldots, M} a_{m}$
		\If{$a_{m^*} > \tau_1$}
			\State add $b_{m^*}$ to the set of estimated change points
			\State NBS$((s, b_{m*}),\{ (\alpha_m,\beta_m)\}_{m=1}^M, \tau_1)$
			\State NBS$((b_{m*}+1,e),\{ (\alpha_m,\beta_m)\}_{m=1}^M,\tau_1)$
			
		\EndIf  
		\OUTPUT The set of estimated change points.
		\caption{Network Binary Segmentation. NBS$((s, e),$ $\{ (\alpha_m,\beta_m)\}_{m=1}^M, \tau_1$)} \label{algorithm:MWBS}
	\end{algorithmic}
\end{algorithm}

\begin{theorem}[Theorem 1 in \citealp{wang2018optimal}]\label{thm-1-network}
Let \Cref{assume:model-network} hold and assume that there exists a constant $C_\alpha>0$ such that, for some $\xi>0$, 
	\[
		\kappa_0 \sqrt{\rho} \ge C_\alpha  \sqrt{\frac{1}{n\Delta}}  \log^{1+ \xi}(T).
	\]
	Let $\{(\alpha_m,\  \beta_m) \}_{m=1}^M\subset (0, T)$ be a collection of intervals whose end points are drawn independently and uniformly from $\{1,\ldots, T\}$ and such that $\max_{m = 1, \ldots, M} (\beta_m -\alpha_m)\le C_R \Delta$, for an absolute constant $C_R>0$.   

Suppose that there exists sufficiently small $0<c_2 <1$ such that the input parameter  $\tau$ of \Cref{algorithm:MWBS} satisfy
	\[
		C_\beta\rho n \log^{3/2}(T)< \tau <c_2\kappa_0^2n^2\rho^2 \Delta. 
	\]
	Then the collection of the estimated change points $\mathcal B=\{\hat \eta_k\}_{k=1}^{\widehat K}$ returned by \Cref{algorithm:MWBS} with input parameters $(0, T)$, $\{(\alpha_m,\beta_m)\}_{m=1}^M$ and $\tau$ is such that	
	\begin{align*}
		& \mathbb{P}\left\{\widehat{K} = K \quad \mbox{and} \quad |\widehat{\eta}_k - \eta_k| \leq C_{\epsilon}\log(T)\left(\frac{\sqrt{\Delta}}{\kappa_k} + \frac{n\rho \log^{1/2}(T)}{\kappa_k^2}\right), \, \forall k \in \{1, \ldots, K\}\right\} \\
		& \hspace{2cm} \ge 1 -\exp\left( \log\frac{T}{\Delta}-M\frac{\Delta^2}{16 T^2} \right) -  T^{-c},
	\end{align*}
	for some absolute constants $C_{\epsilon}, c > 0$ and any $n, T \geq 2$.
\end{theorem}

We remark that, in terms of $\kappa_0$, \Cref{thm-1-network} shows the localisation error is of order
	\[
		\log(T)\left(\frac{\sqrt{\Delta}}{\kappa_0 n\rho} + \frac{\log^{1/2}(T)}{\kappa_0^2 n\rho}\right).
	\]
	Therefore, provided that $\log(T/\Delta) \lesssim M\Delta^2 T^{-2}$, \Cref{thm-1-network} shows that \Cref{algorithm:MWBS} achieves consistent localisation in the regime that
	\[
		\kappa_0 \sqrt{\rho n \Delta} \gtrsim \log^{1+\xi}(T),
	\]
	for any $\xi > 0$.  The role of $\xi$ is the same as that in \Cref{sec-detec-local}.  Together with \Cref{lemma:lower bound testing-network}, we know that it is nearly optimal, save a logarithmic factor.  However, the localisation error achieved in \Cref{thm-1-network} is sub-optimal given \Cref{lem-3.3-lower}.  Two natural questions await: (1) why is it sub-optimal? (2) how can it be improved?
	
The sub-optimality is rooted in the high-dimensionality.  \Cref{algorithm:MWBS} in fact only takes weighted sample mean of matrix inner products, which are merely $\ell_2$-norms of vectorised matrices.  Even for a network of a moderately-high dimension, its vectorised version is of very high dimension.  It is well-understood that merely taking sample means does not lead to good estimation in high-dimensional statistics.  What we learn from \Cref{thm-1-network} is that, if the goal is to localise change points consistently, then one could sacrifice some accuracy in estimating the underlying high-dimensional distributions.  However, if one wishes for more accurate, say optimal change point localisation, then this sacrificed accuracy is probably to be blamed.  

\subsubsection{Optimal localisation}

In order to improve localisation, as we discussed before, one needs to provide more accurate estimation of the underlying distributions.  Just like other problems in high-dimensional statistics, some form of sparsity condition is needed.  In the context of networks, a natural choice of the sparsity is the low rank assumption.  

\begin{assumption}\label{assume:phase 2}
	Let $\{\Theta(t)\}_{t=1}^T$ be defined as in \Cref{assume:model-network}. For some $0 < r \leq n$,
		\[
			\max_{k=1,\ldots,K} \mathrm{rank} \left( \Theta (\eta_{k} ) -\Theta (\eta_{k}-1 ) \right) \le r. 
		\]
\end{assumption}

With the additional low rank assumption \Cref{assume:phase 2}, we will show that the localisation errors can be improved. 

\begin{algorithm}[!ht]
	\begin{algorithmic}
		\INPUT Symmetric matrix $A \in \mathbb{R}^{n \times n}$, $\tau_2, \tau_3 > 0$.
		\State $(\kappa_i(A), v_i) \leftarrow $ the $i$th eigen-pair of $A$, with $|\kappa_1(A)| \geq \cdots |\kappa_n(A)|$
		\State $A' \leftarrow \sum_{i:|\kappa_i (A)|\ge\tau_2} \kappa_i(A) v_iv_i^{\top}$
		\State $\mathrm{USVT}(A, \tau_2,\tau_3) \leftarrow (A''_{ij})$ with
			\[
				(A'')_{ij} \leftarrow \begin{cases}
						(A')_{ij}, & \text{if} \quad |(A'_{ij})| \le \tau_3\\
					\sign ((A')_{ij})\tau_3, & \text{if} \quad |(A'_{ij})| > \tau_3\\
				\end{cases}
			\]
		\OUTPUT $\mathrm{USVT}(A, \tau_2,\tau_3)$.
		\caption{ $\mathrm{USVT}(A, \tau_2,\tau_3)$}
		\label{algorithm:USVT}
	\end{algorithmic}
\end{algorithm} 

\begin{algorithm}[!ht]
	\begin{algorithmic}
		\INPUT $\{X(t)\}_{t=1}^{T}, \{W(t)\}_{t=1}^{T} \in \mathbb{R}^{n\times n}$, $\tau_2, \tau_3$, $\{ \nu_{k}\}_{k=1}^{K} \subset \{2, \ldots, T\}$, $ \nu_0 =0$, $\nu_{K+1} =T+ 1$.		
		\For{$k = 1, \ldots, K$}  
			\State $[s,e] \leftarrow [2^{-1}(\nu_{k-1} + \nu_{k}) , \, 2^{-1}(\nu_{k} + \nu_{k})] $
			\State $\widetilde{\Delta}_k \leftarrow \sqrt {\frac{(e-\nu_k)(\nu_k-s)}{e-s} }  $
			\State $\widehat \Theta_k \leftarrow \mathrm{USVT}(\widetilde B^{s , e} ( \nu_k) , \tau_2,\tau_3\widetilde{\Delta}_k)$
			\State $b_k \leftarrow \argmax_{s\le t \le e  } ( \widetilde A^{ s , e} (t),\widetilde \Theta_k)$
		\EndFor
		\OUTPUT $\{b_k\}_{k=1}^K$.
		\caption{Local Refinement}
		\label{algorithm:RI}
	\end{algorithmic}
\end{algorithm} 

\begin{theorem}[Theorem 2 in \citealp{wang2018optimal}]\label{theorem:localization 1}
Let Assumptions~\ref{assume:model-network} and \ref{assume:phase 2} hold.  Assume that for any $\xi > 0$, there exists an absolute constant $C_\alpha > 0$ such that
	\[
		  \kappa_0 \sqrt{\rho}  \ge C_\alpha \frac{\log^{1 +  \xi}(T)}{\sqrt{\Delta}}\sqrt{\frac{r}{n}}.
	\]
	Let $\mathcal B=\{ \nu_k\}_{k=1}^{ K} \subset \{1,\ldots, T\}$ be a collection of time points.  Suppose that 
	\begin{align}\label{eq-net-cond-init}
		\max_{k=1,\ldots,K} |\nu_k-\eta _k| <\Delta/6.
	\end{align}
	For a large enough absolute constant $C_a > 0$ suppose that 
	\begin{align*}
		\tau_2 = (3/4)(C\sqrt{n\rho} + C_{\varepsilon}\log(T))
		\ \text { and } \
		\tau_3  =\rho,
	\end{align*}    
	where $C > 64 \times 2^{1/4e^2}$ and $C_{\varepsilon} > 12$. Then the outputs of \Cref{algorithm:RI} with input parameters of $(0, T)$, $\{\nu_k\}_{k=1}^K$, $\tau_2$ and $\tau_3$ satisfy
	\begin{align*}
		\mathbb{P}\Bigl\{\max_{k=1,\ldots,K} |\eta_k-\widehat  \eta_k| \le  C_2 \log^2(T) \kappa_0^{-2}n^{-2}\rho^{-1} \Bigr\} \geq 1 - T^{-c},
	\end{align*}
	where $C_2, c > 0$ are absolute constants.
\end{theorem}

\Cref{algorithm:RI} can be seen as a refinement of a set of initial estimators, as we discussed in \Cref{sec-poly-sec-conclusions}.  The refinement is conducted based on a better estimation of the underlying graphons, by using the universal singular value thresholding (USVT, \Cref{algorithm:USVT}) method developed in \cite{Chatterjee2015}.  Note that the output of \Cref{algorithm:MWBS} satisfies the condition on the initial estimators, detailed in \eqref{eq-net-cond-init}, and the signal-to-noise ratio condition required in \Cref{theorem:localization 1} is stronger than that in \Cref{thm-1-network}.  This means \Cref{algorithm:RI} can be used as a second step after \Cref{algorithm:MWBS}, and the final outputs are nearly optimal in terms of the localisation errors.  

Like we discussed in \Cref{sec-poly-sec-conclusions}, one may directly integrate the USVT estimation in the main algorithm \Cref{algorithm:MWBS}, but since the computation costs of conducting the singular value decomposition is of order $O(n^2)$ and the WBS procedure itself is of order $O(n^3)$.  If we directly adopting \Cref{algorithm:USVT} in \Cref{algorithm:MWBS}, then the computational cost is $O(n^5)$.  On the contrary, we use a sample mean in \Cref{algorithm:MWBS} and use the USVT as a refinement, then the computational cost is of order $O(n^4) + O(n^3) = O(n^4)$.

\subsubsection{Conclusions}

The natural question is what happens in the regime 
	\[
		 \sqrt{\frac{1}{n\Delta}}  \log^{1+ \xi}(T) \lesssim \kappa_0 \sqrt{\rho}  \lesssim \sqrt{\frac{r}{n\Delta}}  \log^{1+ \xi}(T).
	\]
	The complete answer is yet known and we will provide some partial answers in line with other high-dimensional statistics problems.

If we replace the subroutine \Cref{algorithm:RI}, which is a polynomial-time algorithm, with an NP-hard graphon-based algorithm \citep[see, e.g.][]{pensky2016dynamic, GaoEtal2015}, then we will be able to produce a nearly optimal localisation rate in the regime 
	\[
		\kappa_0 \sqrt{\rho} \gtrsim  \sqrt{\frac{1 + r^2/n}{n \Delta}}\log^{1+\xi}(T).
	\]
	This means that (i) in the very sparse regime, i.e.~$r \lesssim \sqrt{n}$, the condition required by NP-hard algorithms is nearly optimal, save for a logarithmic factor; (ii) in the moderately sparse regime, i.e.~$\sqrt{n} \lnsim r \lnsim n$, there is a gap between statistical and computational limits; (iii) in the very dense regime, i.e.~$r \asymp n$, NP-hard algorithms are not gaining over polynomial methods.  These observations is consistent with similar phenomena observed in other statistical problems, see e.g.~\cite{zhang2012communication}, \cite{loh2013regularized}, to name but a few. 
	
Finally, the graphon change point analysis has also been studied in different settings over the years, including \cite{zhao2019change}, \cite{LiuEtal2018}, \cite{CribbenYu2017} and \cite{bhattacharjee2018change}, among others.

\subsection{Regression coefficients changes}\label{sec-regression}

\subsubsection{Overview}
In this case, we consider at every time point $t \in \{1, \ldots, T\}$, $(x_t, y_t) \in \mathbb{R}^p \times \mathbb{R}$ is collected, where $y_t$'s are response variables and $x_t$'s are high-dimensional covariates.  In the change point analysis context, we assume the regression coefficients are piecewise constant.  The detailed assumptions are collected below.

\begin{assumption} \label{assume:change point regression model}
Let the data be $\{(x_t, y_t)\}_{t = 1}^T \subset \mathbb{R}^p \times \mathbb{R}$, satisfying
		\[
			y_t = x_t^\top \beta_t^* + \varepsilon_t,
		\]
		where $\beta_t^* \in \mathbb{R}^p$ is the unknown coefficient vector, $x_t$'s are independent and identically distributed, and $\varepsilon_t$'s are independent centred sub-Gaussian random variables with parameters $\sigma_t^2 \leq \sigma_{\varepsilon}^2$ and independent of $\{x_t\}$. 

	In addition, there exists a collection of change points $\{\eta_{k}\}_{k=0}^{K+1} \subset \{1, \ldots, T+1\}$ with $\eta_0 = 1$ and $\eta_{K+1}= T+1$ such that $\beta_t^* \neq \beta_{t - 1}^*$, if and only if $t \in \{\eta_k\}_{k = 1}^K$.
\end{assumption}

\begin{assumption} \label{assume:high dim coefficient}  
	Consider the model defined in \Cref{assume:change point regression model}, where $x_t$'s are centred sub-Gaussian random vectors with $\mathbb{E}(x_t x_t^{\top}) = \Sigma$.  We impose the following additional assumptions.

There exists a subset $S \subset \{1, \ldots, p\}$ such that 
	\[
		\beta_t^*(j) = 0, \quad t = 1, \ldots, T, \quad j \in S^c = \{1, \ldots, p\} \setminus S.
	\]
	Let $d_0 = |S|$.  
			
For some absolute constant $C_{\beta} > 0$, $\max_{t = 1, \ldots, T} \|\beta_t^*\|_{\infty} \leq C_{\beta}$.
	
We have that
	\[
		\Lambda_{\min}(\Sigma) = c_x^2 > 0 \quad \mbox{and} \quad \max_{j = 1, \ldots, p} (\Sigma)_{jj} = C_x^2 > 0.
	\]				
			
Let $\kappa$ and $\Delta$ be the minimal jump size and minimal spacing defined as follows, respectively,
	\[
		\kappa = \min_{k = 1, \ldots, K} \kappa_k = \min_{k = 1, \ldots, K} \| \beta^*_{\eta_k} - \beta^*_{\eta_k - 1}\| \quad \mbox{and} \quad \Delta = \min_{k = 1, \ldots, K + 1} (\eta_{k} - \eta_{k-1}).
	\]
\end{assumption}

The difficulty of this problem is characterised in Lemmas~\ref{lem-reg-lb-1} and \ref{lem-reg-lb-2}, on the minimax lower bounds on detection and localisation, respectively. 

\begin{lemma}[Lemma~3 in \citealp{rinaldo2020localizing}] \label{lem-reg-lb-1}
Let $\{(x_t, y_t)\}_{t = 1}^T \subset \mathbb{R}^p \times \mathbb{R}$ satisfy Assumptions~\ref{assume:change point regression model} and \ref{assume:high dim coefficient}, with $K = 1$.  In addition, assume $x_t \stackrel{\mbox{iid}}{\sim} \mathcal{N}(0, I_p)$ and $\varepsilon_t \stackrel{\mbox{iid}}{\sim} \mathcal{N}(0, \sigma_{\varepsilon}^2)$.  Let $P^T_{\kappa, \Delta, \sigma_{\varepsilon}, d}$ be the corresponding joint distribution.  For any $0 < c < \frac{2}{8e+1}$, consider the class of distributions
	\[
		\mathcal{P}^T = \left\{P^T_{\kappa, \Delta, \sigma_{\varepsilon}, d}: \, \Delta = \min\left\{\lfloor cd_0 \sigma_{\varepsilon}^2\kappa^{-2}\rfloor, \, \lfloor T/4\rfloor \right\}, \, 2cd_0\max\{d_0, 2\} \leq \Delta\right\}.
	\]
	There exists a $T(c)$, which depends on $c$, such that for all $T \geq T(c)$,
	\[
		\inf_{\widehat{\eta}} \sup_{P \in \mathcal{P}^T} \mathbb{E}_P(|\widehat{\eta} - \eta(P)|) \geq \Delta,
	\]
	where $\eta(P)$ is the location of the change point of distribution $P$ and the infimum is over all estimators of the change point. 
\end{lemma}

\begin{lemma}[Lemma~4 in \citealp{rinaldo2020localizing}] \label{lem-reg-lb-2}
Let $\{(x_t, y_t)\}_{t = 1}^T \subset \mathbb{R}^p \times \mathbb{R}$ satisfy Assumptions~\ref{assume:change point regression model} and \ref{assume:high dim coefficient}, with $K = 1$.  In addition, assume $x_t \stackrel{\mbox{iid}}{\sim} \mathcal{N}(0, I_p)$ and $\varepsilon_t \stackrel{\mbox{iid}}{\sim} \mathcal{N}(0, \sigma_{\varepsilon}^2)$.  Let $P^T_{\kappa, \Delta, \sigma_{\varepsilon}, d}$ be the corresponding joint distribution.  For any diverging sequence $\zeta_T$, consider the class of distributions
	\[
		\mathcal{Q}^T = \left\{P^T_{\kappa, \Delta, \sigma, d}: \, \Delta = \min\left\{\lfloor \zeta_T d_0 \sigma_{\varepsilon}^2\kappa^{-2}\rfloor, \, \lfloor T/4\rfloor \right\} \right\}.
	\]
	Then 
	\[
		\inf_{\widehat{\eta}} \sup_{P \in \mathcal{Q}^T} \mathbb{E}_P(|\widehat{\eta} - \eta(P)|)\geq \frac{cd_0 \sigma_{\varepsilon}^2}{\kappa^2},
	\]
	where $\eta(P)$ is the location of the change point of distribution $P$, the infimum is over all estimators of the change point and $c > 0$ is an absolute constant. 
\end{lemma}

Lemmas~\ref{lem-reg-lb-1} and \ref{lem-reg-lb-2} show that in the low signal-to-noise ratio regime $\kappa \sqrt{\Delta} \lesssim \sigma_{\varepsilon} \sqrt{d_0}$, no algorithm is guaranteed to be consistent, and the localisation error lower bound is $d\sigma^2_{\varepsilon} \kappa^{-2}$.

\subsubsection{Consistent localisation}
In order to estimate the change points, we adopt \Cref{algorithm:PDP}.  To be specific, for any interval $I \subset \{1, \ldots, T\}$, let 
	\begin{equation}\label{eq-reg-obj}
		H(I) = \sum_{t \in I}(y_t - x_t^{\top} \widehat \beta^\gamma_I)^2,
	\end{equation}
	with			
	\begin{equation}\label{eq-beta-1}
		\widehat \beta^\gamma_I = \argmin_{v \in \mathbb R^p} \left\{\sum_{t \in I} (y_t - x_t^{\top}v)^2  + \gamma \sqrt{\max\{|I|, \, \log(n \vee p)\}} \| v\|_1\right\},
	\end{equation}
	where $\|\cdot\|_1$ denotes the vector $\ell_1$-norm.  With this construction, the loss function is the sum of residual squares, with a Lasso estimator of the coefficients.  The theoretical guarantees of the output of \Cref{algorithm:PDP} with \eqref{eq-reg-obj} and \eqref{eq-beta-1} are given below.
	
\begin{theorem}[Theorem~1 in \citealp{rinaldo2020localizing}]\label{eq:DP consistent for regression}
Let Assumptions~\ref{assume:change point regression model} and \ref{assume:high dim coefficient} hold.   Assume for any $\xi > 0$, there exists an absolute constant $C_{\mathrm{SNR}} > 0$ such that 
	\begin{equation}\label{eq:snr}
		\Delta \kappa^2 \geq C_{\mathrm{SNR}} d_0^2 K \sigma^2_{\varepsilon}\log^{1 + \xi}(T\vee p).
	\end{equation}
	Let $\{\widetilde{\eta}_k\}_{k = 1}^{\widehat{K}}$ be the output of \Cref{algorithm:PDP}, with the objective function defined in  obtained as solution to the dynamic programming optimisation problem given in \eqref {eq-reg-obj} and \eqref{eq-beta-1} and with tuning parameters
	\[
		\gamma = C_{\gamma}\sigma_{\varepsilon} \sqrt{d_0 \log(T \vee p)} \quad \mbox{and} \quad \lambda = C_{\lambda}\sigma_{\varepsilon}^2 (K+1) d^2_0 \log(T \vee p).
	\]
	It holds that
	\[
		\mathbb{P}\left\{\widehat K = K, \, \max_{k = 1, \ldots, K}|\widetilde{\eta}_k - \eta_k| \leq  \frac{KC_{\epsilon}d_0^2 \sigma^2_{\varepsilon} \log(T \vee p)}{\kappa^2} \right\} \geq 1 - C(T \vee p)^{-c},
	\]
	where $C_{\lambda}, C_{\gamma}, C_{\epsilon}, C, c > 0$ are absolute constants depending only on $C_{\beta}, C_x$ and $c_x$.
\end{theorem}

In view of Lemmas~\ref{lem-reg-lb-1} and \ref{lem-reg-lb-2}, we can see that \Cref{eq:DP consistent for regression} requires a stronger signal-to-noise ratio and achieves a sub-optimal localisation error.  We will improve the localisation in the sequel.  As for the signal-to-noise ratio condition, we remark that if one further assumes $\kappa = O(1)$ and $K = O(1)$, then one can replace \eqref{eq:snr} with
	\[
		\Delta \kappa^2 \geq C_{\mathrm{SNR}} d_0 \sigma^2_{\varepsilon}\log^{1 + \xi}(T\vee p),
	\]
	and \Cref{eq:DP consistent for regression} still holds.  This shows the nearly-optimality of \Cref{algorithm:PDP} in terms of the signal-to-noise ratio condition.  However, without the extra condition that $\kappa = O(1)$ and $K = O(1)$, it remains an \textbf{open problem} in deriving a consistent change point estimator under minimal conditions.
	
\subsubsection{Optimal localisation}

So far we have already used the refinement idea twice.  In \Cref{sec-piecewise-polynomial}, a refinement step is to improve the localisation rate so it is not a function of $K$, the number of change points.  In \Cref{sec-graphon}, a refinement step is to improve the localisation rate by providing a better estimation of the underlying high-dimensional objects.  Comparing the localisation error in \Cref{eq:DP consistent for regression} and the minimax lower bound in \Cref{lem-reg-lb-2}, we see that a refinement should ideally eliminate the dependence on $K$ and improve from $d_0^2$ to $d_0$.  This suggests that a refinement should not only work in the intervals containing one and only one true change point, but also need to provide better estimation of the underlying distributions.  
	
\begin{algorithm}[htbp]
\begin{algorithmic}
	\INPUT Data $\{(x_t, y_t)\}_{t=1}^T$, a collection of time points $\{\widetilde{\eta}_k\}_{k = 1}^{\widetilde{K}}$ , tuning parameter $\zeta > 0$.
	\State $(\widetilde{\eta}_0, \widetilde{\eta}_{\widetilde{K} + 1}) \leftarrow (0, T)$
	\For{$k = 1, \ldots, \widetilde{K}$}  
		\State $(s_k, e_k) \leftarrow (\widetilde{\eta}_{k-1}/3 + 2\widetilde{\eta}_{k}/3, 2\widetilde{\eta}_{k}/3 + \widetilde{\eta}_{k+1}/3)$
		\State 
		\begin{align}
			\left(\widehat{\beta}_1, \widehat{\beta}_2, \widehat{\eta}_k\right) \leftarrow \argmin_{\substack{\eta \in \{s_k + 1, \ldots, e_k - 1\} \\ \beta_1, \beta_2 \in \mathbb{R}^{p} \\ \beta_1 \neq \beta_2}}   \Bigg\{\sum_{t = s_k + 1}^{\eta}\bigl\|y_{t} - \beta_1^{\top} x_t\bigr\|^2 + \sum_{t = \eta + 1}^{e_k}\bigl\|y_t - \beta_2 x_t\bigr\|^2 \nonumber \\
		+ \zeta  \sum_{i = 1}^p \sqrt{(\eta - s_k)(\beta_1)_{i}^2 + (e_k - \eta)(\beta_2)_{i}^2}\Bigg\} \label{eq-linear-refine-criterion}
		\end{align}
	\EndFor
	\OUTPUT The set of estimated change points $\{\widehat{\eta}_k\}_{k = 1}^{\widetilde{K}}$.
\caption{Local refinement. LR$(\{(x_t, y_t)\}_{t=1}^T, \{\widetilde{\eta}_k\}_{k = 1}^{\widetilde{K}}, \zeta)$}
\label{algorithm:LR-regression}
\end{algorithmic}
\end{algorithm} 
	
\begin{theorem}[Corollary 2 in \citealp{rinaldo2020localizing}] \label{cor-lr-high-dim}
Assume the same conditions of \Cref{eq:DP consistent for regression}.  Let $\{\widetilde{\eta}_k\}_{k = 1}^{K}$ be a set of time points satisfying
	\begin{equation}\label{eq-lr-cond-1-linear}
		\max_{k = 1, \ldots, K} |\widetilde{\eta}_k - \eta_k| \leq \Delta/7.
	\end{equation}
	Let $\{\widehat{\eta}_k\}_{k = 1}^{\widehat{K}}$ be the change point estimators generated from \Cref{algorithm:LR-regression} with $\{\widetilde{\eta}_k\}_{k = 1}^{K}$ and 
	\[
		\zeta = C_{\zeta}\sqrt{\log(T \vee p)}
	\]
	as inputs.  Then,
	\[
		\mathbb{P}\left\{\widehat K = K, \, \max_{k = 1, \ldots, K}|\hat{\eta}_k - \eta_k| \leq \frac{C_{\epsilon} d_0 \log(T \vee p)}{\kappa^2} \right\} \geq 1 - T^{-c},
	\]
	where $C_{\zeta}, C_{\epsilon}, c > 0$ are absolute constants depending only on $C_{\beta}, \mathcal{M}$ and $c_x$.
\end{theorem}

\Cref{cor-lr-high-dim} shows that if \Cref{algorithm:LR-regression} is adopted as a refinement step of \Cref{algorithm:PDP} with \eqref{eq-reg-obj} and \eqref{eq-beta-1}, then the corresponding localisation error is nearly-optimal, off by a logarithmic factor.  The near optimality is achieved due to two key ingredients.
	\begin{itemize}
	\item The condition on the initial estimators \eqref{eq-lr-cond-1-linear} provides the opportunity that we are able to work in intervals containing one and one true change point.  This eliminates the dependence on $K$ in the localisation error.
	\item In \Cref{algorithm:LR-regression}, a group lasso estimation is adopted in \eqref{eq-linear-refine-criterion}.  Since we have already guaranteed that there is only one change point in the working interval, the group lasso penalty captures this feature and returns better estimation, with a higher computational cost.
	\end{itemize}

In this section, we only reviewed the coefficients change in linear regression models.  In fact, similar techniques can also be extended to other forms of regression problems, including (vector) autoregressive models, self-exciting Poisson processes, and other time series models.   We conclude this subsection with a list of existing literature on different aspects of different regression coefficient change point problems.  These papers include \cite{aue2006change}, \cite{wang2019statistically}, \cite{wang2020detecting}, \cite{safikhani2020joint}, \cite{leonardi2016computationally} and others.

\subsection{Conclusions}	\label{sec-high-dim-conclusion}

The three high-dimensional cases we reviewed here are representative.  
	\begin{itemize}
	\item In \Cref{sec-cov-change}, despite the high-dimensionality, we are able to find a polynomial-time algorithm achieves nearly optimal localisation rates under nearly optimal signal-to-noise conditions, both off by logarithmic factors.  In other words, all three goals we listed in \Cref{sec-what-we-will-cover} are achieved.
	\item In \Cref{sec-graphon}, we reviewed a case exhibiting statistical and computational tradeoffs.  In terms of the three goals we listed in \Cref{sec-what-we-will-cover}, only the third one is achieved.
	\item In \Cref{sec-regression}, we show that under some mild conditions and a nearly optimal signal-to-noise ratio condition, a penalisation-based method is able to provide consistent change point estimators, based on which, a refinement can improve the localisation error to be nearly optimal.  In terms of the three goals we listed in \Cref{sec-what-we-will-cover}, only the third one is achieved.  With some mild extra conditions, all three goals can be achieved.
	\end{itemize}

Recall that in the univariate mean change point analysis in \Cref{sec-detec-local}, there exist nearly-optimal polynomial-time methods, in terms of both detection and localisation.  The high-dimensionality obviously increases difficulties.  The reason that we can achieve the near optimality without any additional steps or conditions in the covariance change point problem in \Cref{sec-cov-change}, is largely due to the matrix operator norm used in the model assumption.  The operator norm plays the role of dimension reduction -- it essentially means all the useful information lies in the largest eigenvalue, despite the high-dimensionality of the data.  To elaborate, if instead of the operator norm, we use the entry-wise maximum norm, then we can still achieve the near-optimality despite the high-dimensionality.  However, if we use the Frobenius norm instead of the operator norm, then we will essentially meet the difficulty as that in \Cref{sec-graphon}.

In \Cref{sec-graphon}, the matrix Frobenius norm is adopted in defining the distributional differences.  Different from the matrix operator norm, the Frobenius norm is not helping at all in terms of dimension reduction.  In other words, the jump size $\kappa$ defined thereof is allowed to vary in $(0, n\rho)$.  The high-dimensional nature of the problem leads to the gaps in \Cref{sec-graphon}.

We conclude this section with a list of papers working on change point detection in other types of high-dimensional data.  
\begin{itemize}
\item High-dimensional mean change points: \cite{AstonKirch2014}, \cite{BarigozziEtal2016}, \cite{Cho2015}, \cite{ChoFryzlewicz2015}, \cite{HorvathHuskova2012}, \cite{wang2016high} and \cite{Jirak2015}, among others.
\item Graphical models and high-dimensional precision matrices change points: \cite{amini2013sequential}, \cite{gibberd2017regularized}, \cite{gibberd2014high}, \cite{gibberd2017multiple}, \cite{keshavarz2018sequential} and \cite{londschien2019change}, among others.
\item Functional data change point analysis:  \cite{BerkesEtal2009}, \cite{aue2009estimation}, \cite{aston2012detecting}, \cite{li2018bayesian}, \cite{aston2011power}, \cite{chiou2019identifying}, \cite{jiao2020break}, \cite{liu2020fast} and \cite{dette2019detecting}, among others.
\end{itemize}

\section{Extension 3: Nonparametric models}\label{sec-para-nonpar}

In this section, we discuss nonparametric models.  We will build up from a univariate case, then a multivariate case and conclude with a reproducing kernel Hilbert space case.  Regarding the distances used to characterise the distribution differences, we will cover three different distances.  As for the algorithms, we will study both the penalisation-based and scan-statistics-based methods.

\subsection{Univariate}\label{sec-np-uni}

Different from the cases studied in Sections~\ref{sec-detec-local} and \ref{sec-piecewise-polynomial}, in this subsection, the differences of the distributions are characterised by the Kolmogorov--Smirnov distance.  

\begin{assumption}\label{assump-uni-np}
Let $\{X_t\}_{t = 1}^T \subset \mathbb{R}$ be a collection of independent random variables such that $X_t \sim F_t$, where $F_t$'s are cumulative distribution functions (CDFs).  Let $\{\eta_k\}_{k = 0}^{K+1} \subset \{1, \ldots, T+1\}$ be a collection of change points with $1 = \eta_0 < \eta_1 < \ldots < \eta_K \leq T < \eta_{K+1} = T+1$ such that
	\[
		F_t \neq F_{t - 1}, \quad \mbox{if and only if } t \in \{\eta_1, \ldots, \eta_K\}.
	\]

Define the minimal spacing $\Delta$ and the jump size $\kappa$ as	
	\[
		\Delta = \min_{k = 1, \ldots, K+1} \{\eta_k - \eta_{k-1}\} > 0
	\]
	and
	\begin{equation}\label{eq-as1-kappa}
		\kappa = \min_{k = 1, \ldots, K} \kappa_k = \min_{k = 1, \ldots, K} \sup_{z \in \mathbb{R}} \bigl|F_{\eta_{k}}(z) - F_{\eta_{k}-1}(z)\bigr| > 0.
	\end{equation}	
\end{assumption}

The difficulty of this problem is characterised in Lemmas~\ref{lemma-low-snr-uninp} and \ref{lemma-error-opt-uni-np}, focusing on the minimax lower bounds on detection and localisation, respectively.

\begin{lemma}[Lemma 3 in \citealp{padilla2019optimal}]\label{lemma-low-snr-uninp}
Let $\{X_t\}_{t = 1}^T$ be a time series satisfying \Cref{assump-uni-np} with one and only one change point.  Let $P^T_{\kappa, \Delta}$ denote the corresponding joint distribution.  For any $0 < \zeta < 1/\sqrt{2}$, denote
	\[
		\mathcal{P}^T = \left\{P^T_{\kappa, \Delta}: \, \Delta = \min\left\{\left\lfloor \frac{\zeta^2}{\kappa^2} \right\rfloor, \, \left\lfloor \frac{T}{3} \right\rfloor \right\}\right\}.
	\]
	Let $\widehat{\eta}$ and $\eta(P)$ be an estimator and the true change point, respectively.  It holds that
	\[
		\inf_{\widehat{\eta}} \sup_{P \in \mathcal{P}^T} \mathbb{E}_P\bigl(\bigl|\hat{\eta} - \eta(P)\bigr|\bigr) \geq (1 - 2\zeta^2)\Delta,
	\]
	where the infimum is over all possible estimators of the change point location.
\end{lemma}

\begin{lemma}[Lemma 4 in \citealp{padilla2019optimal}]\label{lemma-error-opt-uni-np}
Let $\{X_t\}_{t = 1}^T$ be a time series satisfying \Cref{assump-uni-np} with one and only one change point.  Let $P^T_{\kappa, \Delta}$ denote the corresponding joint distribution.  Consider the class of distributions
	\[
		\mathcal{Q}^T = \left\{P^T_{\kappa, \Delta}: \, \Delta < T/2,\, \kappa < 1/2, \, \kappa\sqrt{\Delta} \geq \zeta_T \right\},
	\]
	for any sequence $\{ \zeta_T \}$ such that $\lim_{T \rightarrow \infty} \zeta_T = \infty $. Let $\widehat{\eta}$ and $\eta(P)$ be an estimator and the true change point, respectively.  Then, for all $T$ large enough, it holds that 
	\[
		\inf_{\widehat{\eta}} \sup_{P \in \mathcal{Q}^T} \mathbb{E}_P\bigl(\bigl|\widehat{\eta} - \eta(P)\bigr|\bigr) \geq \max \left\{ 1, \frac{1}{2} \Big\lceil\frac{1}{\kappa^2} \Big\rceil e^{-2} \right\},
	\]
	where the infimum is over all possible estimators of the change point locations.	
\end{lemma}

Lemmas~\ref{lemma-low-snr-uninp} and \ref{lemma-error-opt-uni-np} show that in the low signal-to-noise ratio regime $\kappa\sqrt{\Delta} \lesssim 1$, no algorithm is guaranteed to provide consistent change point estimators, and the minimax lower bounds on the localisation is $\kappa^{-2}$.  

We will demonstrate how a scan-statistics-based method is able to reach near optimality in the sense of both detection and localisation.   Based on \Cref{assump-uni-np}, we tailor the CUSUM statistics defined in \Cref{def-cusum} to incorporate the Kolmogorov--Smirnov distance. \Cref{def:cusum-KS} in fact replaces $X_t$'s in \Cref{def-cusum} with indicator functions $\mathbbm{1}\{X_t \leq s\}$, $t \in \{1, \ldots, T\}$, $s \in \mathbb{R}$.

\begin{definition} \label{def:cusum-KS}
	For any integer triplet $(s, t, e)$, $0 \leq s < t < e \leq T$, define 
	\[
		D_{s,e}^t = \sup_{z \in \mathbb{R}} \left|\sqrt{\frac{(t-s) (e-t)}{e-s}}\left\{\widehat{F}_{(s+1): t}(z) - \widehat{F}_{(t+1):e}(z)\right\}\right|,
	\]
	where for all integer pair $(s, e)$, $1 \leq s < e \leq T$ and any $z\in  \mathbb{R}$,
	\[
		\widehat{F}_{s: e}(z) = \frac{1}{e-s} \sum_{t=s}^{e}\mathbbm{1}_{\{X_t \leq z\}}.
	\]
\end{definition}

With the Kolmogorov--Smirnov version of the CUSUM statistics, we can adapt \Cref{algorithm:WBS} by replacing the CUSUM statistics there with the one in \Cref{def:cusum-KS}.  To be specific, given data $\{X_t\}_{t = 1}^T$ and any integer triplet $(s, t, e)$, $0 \leq s < t < e \leq T$, we let
	\begin{equation}\label{eq-dtilde-uni-np}
		\widetilde{D}^t_{s, e} = \max_{i = 1, \ldots, T} \left|\sqrt{\frac{(t-s) (e-t)}{e-s}}\left\{\widehat{F}_{(s+1): t}(X_i) - \widehat{F}_{(t+1):e}(X_i)\right\}\right|,
	\end{equation}
	i.e.~the supreme is taken on the support of all data points.  The theoretical guarantee is give below.
	
\begin{theorem}[Theorem 2 in \citealp{padilla2019optimal}]\label{thm-wbs-np-uni}
Let the CUSUM statistics used in \Cref{algorithm:WBS} be \eqref{eq-dtilde-uni-np}.  Assume the inputs of \Cref{algorithm:WBS} are as follows:
	\begin{itemize}
	\item the sequence $\{X_t\}_{t = 1}^T$ satisfies \Cref{assump-uni-np} and there exists a constant $C_{\mathrm{SNR}} > 0$ such that
		\[
			\kappa \sqrt{\Delta} > C_{\mathrm{SNR}} \sqrt{\log(T)};
		\]	
	\item the collection of intervals $\{(\alpha_m, \beta_m)\}_{m = 1}^M \subset \{1, \ldots, T\}$, with endpoints drawn independently and uniformly from $\{1, \ldots, T\}$, satisfy $\max_{m = 1, \ldots, M}(\beta_m - \alpha_m) \leq C_R\Delta$, almost surely, for an absolute constant $C_R > 1$; and 
	\item the tuning parameter $\tau$ satisfies $c_{\tau, 1}\sqrt{\log(T)} \leq \tau \leq c_{\tau, 2}  \kappa \Delta^{1/2}$,	 where $c_{\tau, 1}, c_{\tau, 2} > 0$ are constants.
	\end{itemize}

Let $\{\widehat{\eta}_k\}_{k = 1}^{\widehat{K}}$ be the corresponding output.  Then
	\begin{align}
		& \mathbb{P}\left\{\widehat{K} = K \quad \mbox{and} \quad  |\widehat{\eta}_k - \eta_k| = \epsilon_k \leq C_{\epsilon}\kappa_k^{-2}\log(T), \, \forall k = 1, \ldots, K\right\} \nonumber\\
		& \hspace{1cm} \geq 1 - 24 \log(T) T^{-4} - \frac{48 }{\log(T) \Delta} - \exp\left\{\log\left(\frac{T}{\Delta}\right) - \frac{M\Delta^2}{16T^2}\right\}, \nonumber
	\end{align}
	where $C_{\epsilon} > 0$ is an absolute constant.
\end{theorem}	

\Cref{thm-wbs-np-uni} shows that there exists a CUSUM-based algorithm which is nearly optimal in terms of both detection and localisation, save for logarithmic factors.  

It is interesting to compare the results we obtained here and those in \Cref{sec-detec-local}.  Comparing the signal-to-noise ratio conditions on consistent localisations
	\[
		\kappa \sqrt{\Delta} \gtrsim \log(T) \quad \mbox{and} \quad \kappa \sqrt{\Delta} \gtrsim \sigma \log(T),
	\]
	and the two localisation errors
	\[
		\kappa^{-2} \quad \mbox{and} \quad \kappa^{-2} \sigma^2,
	\]
	we see that the univariate nonparametric change point detection problem with Kolmogorov--Smirnov distance, can be seen as a univariate mean change point detection problem with $\sigma \asymp O(1)$.  This connection is due to the definition of empirical distribution functions used in \Cref{def:cusum-KS}.  Each observation is turned to an indicator variable, which is a Bernoulli random variable with variance upper bounded by 1.
	
Regarding the comparisons with \Cref{sec-detec-local}, another remark is in \Cref{sec-detec-local}, it is known that a minimax lower bound on detection is
	\[
		\kappa\sqrt{\Delta} \lesssim \sigma \sqrt{\log(T)},
	\]
	while in this subsection, the state-of-the-art result on the lower bound is
	\[
		\kappa\sqrt{\Delta} \lesssim 1,
	\]
	which leaves a gap of order $\log^{1/2+\xi}(T)$ between the lower and upper bounds.  It remains an \textbf{open question} on how to close this gap even further.  We conjecture that the gap is due to a loose lower bound.

\subsection{Multivariate}\label{sec-np-multi}

In the second nonparametric example, we study a sequence of random vectors and the distance used to define the distribution jumps is chosen to be the supreme norm of a function.  The detailed model assumption is provided below.

\begin{assumption} \label{assump-model-multi-np}

Let $\{X_t\}_{t = 1}^T \subset \mathbb{R}^p$ be a sequence of independent random vectors with unknown distributions $\{P_t \}_{t = 1}^T$ such that, for an unknown sequence of change points $\{ \eta_k\}_{k=1}^K \subset \{2, \ldots, T\}$ with $1 = \eta_0 < \eta_1 < \ldots < \eta_K \leq T < \eta_{K+1} = T+1$, we have 
	\[
		P_{t} \neq P_{t-1} \quad \mbox{if and only if} \quad t \in \{\eta_1, \ldots, \eta_K \}.
	\]

Assume that, for each $t =1 ,\ldots, T$, the distribution $P_t$ has a bounded Lebesgue  density function $f_t: \, \mathbb{R}^p \to \mathbb{R}$ such that 
	\[
		\max_{t = 1, \ldots, T} \bigl|f_t(s_1) - f_t(s_2)\bigr| \leq C_{\mathrm{Lip}} \|s_1 - s_2\|, \quad \mbox{for all } s_1, s_2 \in \mathcal{X},
	\]
	where  $\mathcal{X} \subset \mathbb{R}^p$ is the union of the supports of all the density functions $f_t$, $\|\cdot\|$ represents the $\ell_2$-norm, and $C_{\mathrm{Lip}} > 0$ is an absolute constant. 	
	We let 
	\[
		\Delta = \min_{k = 1, \ldots, K+1} \{\eta_k - \eta_{k-1}\}
	\]
	denote the minimal spacing between any two consecutive change points, and let
	\[
		\kappa = \min_{k = 1, \ldots, K} \kappa_k = \min_{k = 1, \ldots, K} \sup_{z \in \mathbb{R}^p} \bigl|f_{\eta_{k}}(z) - f_{\eta_{k}-1}(z)\bigr| =  \|f_{\eta_k} - f_{\eta_{k}-1}\|_{\infty} > 0
	\]
	be the minimal jump size.
\end{assumption}

As usual, we first study the minimax lower bounds determining the difficulties of the problem.

\begin{lemma}[Lemma 2 in \citealp{padilla2019optimal2}]\label{lem-snr-lb-multi-np}
Let $\{X_t\}_{t = 1}^{T}$ be a sequence of random vectors satisfying Assumption~\ref{assump-model-multi-np} with one and only one change point and let  $P^T_{\kappa,  \Delta}$ denote the corresponding joint distribution. Then, there exist universal positive constants $C_1$, $C_2$ and $c < \log(2)$ such that, for all $T$ large enough,
	\[
		\inf_{\widehat{\eta}} \sup_{P \in \mathcal{P}^T} \mathbb{E}_P\bigl(\bigl|\widehat{\eta} - \eta(P)\bigr|\bigr) \geq \Delta/4,
	\]
	where 
	\[
	\mathcal{P}^T = \left\{P^T_{\kappa, \Delta}: \, \Delta < T/2,\,\kappa < C_1, \, \kappa^{p+2}  \Delta \leq c, \, C_{\mathrm{Lip}}\leq  C_2\right\},
	\]
	the quantity $\eta(P)$ denotes the true change point location of $P \in \mathcal{P}^T$ and the infimum is over all possible estimators of the change point location.
\end{lemma}

\begin{lemma}[Lemma 3 in \citealp{padilla2019optimal2}]\label{lemma-error-opt-multi-np}
	Let $\{X_t\}_{t = 1}^{T}$ be a sequence of random vectors satisfying Assumption~\ref{assump-model-multi-np} with one and only one change point and let  $P^T_{\kappa,  \Delta}$ denote the corresponding joint distribution. Then, there exist universal positive constants $C_1$ and $C_2$ such that, for any sequence $\{ \zeta_T \}$ satisfying $\lim_{T \rightarrow \infty} \zeta_T = \infty $, 
\[
	\inf_{\hat{\eta}} \sup_{P \in \mathcal{Q}} \mathbb{E}_P\bigl(\bigl|\hat{\eta} - \eta(P)\bigr|\bigr) \geq \max \left\{ 1, \frac{1}{4} \Big\lceil\frac{1}{V_p^2\kappa^{p+2} } \Big\rceil e^{-2} \right\},
	\]
where $V_p = \pi^{p/2} (\Gamma(p/2 + 1))^{-1}$ is the volume of a unit ball in $\mathbb{R}^p$, 
	\[
	\mathcal{Q}^T = \left\{P^T_{\kappa, \Delta}: \, \Delta < T/2,\,\kappa < C_1, \, \kappa^{p+2}  V_p^2 \Delta \geq \zeta_T, \, C_{\mathrm{Lip}}\leq  C_2\right\},
	\]
	the quantity $\eta(P)$ denotes the true change point location of $P \in \mathcal{Q}^T$ and the infimum is over all possible estimators of the change point location.
\end{lemma}

Lemmas~\ref{lem-snr-lb-multi-np} and \ref{lemma-error-opt-multi-np} show that in the low signal-to-noise ratio regime where $\kappa^{p+2}\Delta \lesssim 1$, no algorithm is guaranteed to be consistent and a minimax lower bound on the localisation error is of order $\kappa^{-(p+2)}$.

In order to match these lower bound, we use the CUSUM-based methods again here and summon \Cref{algorithm:WBS} by adjusting the CUSUM statistics.  Following the same routine, we first define the corresponding CUSUM statistics.  \Cref{def-mul-non-cusum} replaces the data $X_t$'s in \Cref{def-cusum} with a kernel function. 

\begin{definition}\label{def-mul-non-cusum}
	Let $\{X_t\}_{t=1}^T$ be a sample in $\mathbb{R}^p$. For any integer triplet $(s, t, e)$ satisfying $0 \leq s < t < e \leq T$ and any $x \in \mathbb{R}^p$,  the multivariate nonparametric CUSUM statistic is defined as the function
	\[
	x \in \mathbb{R}^p \mapsto	\widetilde{Y}^{s, e}_{t}(x) = \sqrt{\frac{(t-s)(e-t)}{e-s}} \left\{\hat{f}_{s+1, t, h}(x) - \hat{f}_{t+1, e, h}(x)\right\},
	\]
	where
	\[
		\hat{f}_{s, e, h}(x) = \frac{h^{-p}}{e - s} \sum_{i = s + 1}^e \mathpzc{k} \left(\frac{x - X(i)}{h}\right)
	\]
	and $\mathpzc{k}(\cdot)$ is a kernel function \citep[see e.g.][]{parzen1962estimation}. In addition, define	
	\[
		\widetilde{Y}^{s, e}_t = \max_{i = 1, \ldots, T}\left|\widetilde{Y}^{s, e}_t (X_i)\right|.
	\]
\end{definition}

Note that in \Cref{def-mul-non-cusum}, the CUSUM statistics is based on a kernel estimator of underlying densities.  The theoretical guarantees of \Cref{algorithm:WBS} with \Cref{def-mul-non-cusum} are presented in \Cref{thm-wbs-np-multi}, with additional assumptions collected in \Cref{assump-kernel}.

\begin{assumption}\label{assump-kernel}
	 Let $\mathpzc{k}: \, \mathbb{R}^p \to \mathbb{R}$ be a kernel function with $\|\mathpzc{k}\|_{\infty}, \|\mathpzc{k}\|_2 < \infty$ such that,
	\begin{itemize}
	\item [(i)] the class of functions
		\[
			\mathcal{F}_{\mathpzc{k}, [l, \infty)} = \left\{\mathpzc{k}\left(\frac{x - \cdot}{h}\right): \, x \in \mathcal{X}, h \geq l\right\}
		\]
		from $\mathbb{R}^p$ to $\mathbb{R}$ is separable in $L_{\infty}(\mathbb{R}^p)$, and is a uniformly bounded VC-class with dimension $\nu$, i.e.~there exist positive numbers $A$ and $\nu$ such that, for every positive measure $Q$ on $\mathbb{R}^p$ and for every $u \in (0, \|\mathpzc{k}\|_{\infty})$, it holds that
		\[
			\mathcal{N}(\mathcal{F}_{\mathpzc{k}, [l, \infty)}, L_2(Q), u) \leq \left(\frac{A\|\mathpzc{k}\|_{\infty}}{u}\right)^{\nu};
		\]
	\item [(ii)] for a fixed $m > 0$, 
		\[
			\int_0^{\infty} t^{p-1} \sup_{\|x\| \geq t} |\mathpzc{k}(x)|^m\, dt < \infty.
		\]
	\item [(iii)] there exists a constant $C_{\mathpzc{k}} > 0$ such that
		\[
			\int_{\mathbb{R}^p}\mathpzc{k}(z) \|z\| \, dz \leq C_{\mathpzc{k}}.
		\]
	\end{itemize}
\end{assumption}

\begin{theorem}[Theorem 1 in \citealp{padilla2019optimal2}] \label{thm-wbs-np-multi}
Assume that the sequence $\{X_t \}_{t=1}^T$ satisfies the model described in \Cref{assump-model-multi-np} and assume that for a given $\xi > 0$, there exists an absolute constant $C_{\mathrm{SNR}} > 0$ such that
	\[
		\kappa^{p+2}  \Delta > C_{\mathrm{SNR}} \log^{1 + \xi}(T).  
	\]

Let $\mathpzc{k}(\cdot)$ be a kernel function satisfying \Cref{assump-kernel}. 
Then, there exist positive universal constants $C_R$, $c_{\tau, 1}$, $c_{\tau, 2}$ and $c_h$, such that if \Cref{algorithm:WBS} is applied to the sequence $\{X_t \}_{t=1}^T$ using the CUSUM statistics defined in \Cref{def-mul-non-cusum}, any collection  $\{(\alpha_m, \beta_m)\}_{m = 1}^M \subset \{1, \ldots, T\}$ of random time intervals with endpoints drawn independently and uniformly from $\{1, \ldots, T\}$ with $\max_{m = 1, \ldots, M}(\beta_m - \alpha_m) \leq C_R\Delta$,  almost surely, tuning parameter $\tau$ satisfying 
	\[
		c_{\tau, 1}\max\left\{h^{-p/2}\log^{1/2}(T), \, h \Delta^{1/2} \right\} \leq \tau \leq c_{\tau, 2}  \kappa \Delta^{1/2},
	\] 
	and bandwidth $h$ given by $h = c_{h} \kappa$, then the resulting change point estimator $\{ \widehat{\eta}_k\}_{k=1}^{\widehat{K}}$  satisfies
		\begin{align}
		& \mathbb{P}\left\{\widehat{K} = K \quad \mbox{and} \quad \epsilon_k = |\hat{\eta}_k - \eta_k| \leq C_{\epsilon}\kappa^{-2}_k \kappa^{-p} \log(T), \, \forall k = 1, \ldots, K\right\} \nonumber\\
		& \hspace{2cm} \geq  1 -  3T^{-c} - \exp\left\{\log\left(\frac{T}{\Delta}\right) - \frac{M\Delta}{4C_RT}\right\}, \nonumber
		\end{align}
		for universal positive constants $C_\epsilon$ and $c$.
\end{theorem}

\Cref{thm-wbs-np-multi} shows that a CUSUM-based method is nearly optimal in terms of both detection and localisation.  

Recall in \Cref{sec-graphon} we remark that if the goal is to estimate change points, then we can sacrifice some accuracy in estimating the underlying distributions.  We actually have similar observations here.  In the density estimation literature, with the Lipschitz condition imposed on the densities, the optimal bandwidth rate is $h_1 \asymp \{\log(\Delta)/\Delta\}^{1/(p+2)}$, with the goal of estimating the underlying densities.  In \Cref{thm-wbs-np-multi}, we see that in order to obtain optimal change point estimation, the bandwidth is required to be $h_{\mathrm{opt}} \asymp \kappa$.  In the situations where $\kappa \gtrsim \{\log(\Delta)/\Delta\}^{1/(p+2)}$, using a bandwidth $h_{\mathrm{opt}}$ lead to larger bias in estimating the densities, but is required if the goal is to estimate change points optimally.

We reviewed the state-of-the-art results in this subsection.  The dimensionality $p$ is considered as an absolute constant.  It remains an \textbf{open problem} if we allow $p$ to diverge, what the minimax rates are in terms of both detection and localisation.  

\subsection{A reproducing kernel Hilbert space}\label{sec-general-np}

So far the random objects we considered are in Euclidean spaces and the optimal methods we present are both scan-statistics-based methods.  In this section, we consider general $\mathcal{X}$-valued random objects, where $\mathcal{X}$ is an arbitrary (measurable) space, and consider a penalisation-based-method.  The change points are defined to be the change points in a reproducing kernel Hilbert space, which is induced by a certain kernel.  The detailed assumptions are collected below.

\begin{assumption}\label{assump-general-np}
Let $\{X_t\}_{t = 1}^T \subset \mathcal{X}$ be a sequence of independent random objects with unknown distributions $\{P_t\}_{t = 1}^T$, where $\mathcal{X}$ is a measurable space.

Let $\mathpzc{k}: \, \mathcal{X} \times \mathcal{X} \to \mathbb{R}$ be a positive semidefinite kernel.  Let $\mathcal{H}$ be the reproducing kernel Hilbert space associated with the kernel $\mathpzc{k}$, together with the canonical feature map $\Phi: \, \mathcal{X} \to \mathcal{H}$, satisfying $\Phi(x) = \mathpzc{k}(\cdot, x)$, $x \in \mathcal{X}$.   Assume $\mathcal{H}$ is separable.  For any $t \in \{1, \ldots, T\}$, define $Y_t = \Phi(X_t) \in \mathcal{H}$, let $\mu_t^*$ be the Bochner integral \citep[e.g.][]{ganiev2013bochner} of $Y_t$ and $\varepsilon_t = Y_t - \mu_t^*$.  Assume that there exists a positive absolute constant $V$ such that
	\[
		\max_{t = 1, \ldots, T} \mathbb{E}\left(\|\varepsilon_t\|_{\mathcal{H}}^2\right) \leq V.
	\]	
	
Let $\{\eta_k\}_{k = 0}^{K+1} \subset \{1, \ldots, T+1\}$ be a strictly increasing sequence, with $\eta_0 = 1$ and $\eta_{K+1} = T+1$, satisfying
	\[
		\mu_t^* \neq \mu_{t-1}^* \quad \mbox{if and only if } t \in \{\eta_k\}_{k = 1}^K. 
	\]
	Let
	\[
		\kappa = \min_{k = 1, \ldots, K}\|\mu^*_{\eta_k} - \mu^*_{\eta_k - 1}\|_{\mathcal{H}} \quad \mbox{and} \quad \Delta = \min_{k = 1, \ldots, K + 1} (\eta_k - \eta_{k-1}),
	\]	
	where $\|\cdot\|_{\mathcal{H}}$ is the norm of $\mathcal{H}$.
\end{assumption}

In order to handle random objects in general space, \Cref{assump-general-np} adopts a kernel function to turn the general space into univariate random variables.  Recall in \Cref{sec-intro}, the change points are defined to be the change points of the sequence $\{P_t\}$, but the change points in \Cref{assump-general-np} are defined to be the change points of $\{\mu_t^*\}$.  Note that, if the kernel $\mathpzc{k}(\cdot, \cdot)$ is a characteristic kernel, then $X_t$ and $X_{t+1}$ have the same distribution if and only if $\mu_t^* = \mu_{t+1}^*$.  In this case, the set of the change points of $\{\mu_t^*\}$ is identical to that of the change points of $\{P_t\}$.  In general, this one-to-one correspondence might not hold.  

In order to understand the difficulties of this problem, we see that the univariate mean change point problem we studied in \Cref{sec-detec-local} is a special case of \Cref{assump-general-np} with $\mathcal{X} = \mathbb{R}$, $\mathpzc{k}(x, y) = xy$, $x, y \in \mathbb{R}$, and assume $\sigma = O(1)$.  This shows that the minimax lower bounds we achieved in \Cref{sec-detec-local} are still valid here: in the low signal-to-noise ratio regime $\kappa \sqrt{\Delta} \gtrsim \log^{1/2}(T)$, no algorithm is guaranteed to be consistent, and a minimax lower bound on the localisation error is $\kappa^{-2}$.

To match the lower bounds, we resort to the penalisation-based method in \Cref{algorithm:PDP} with the function $H(I)$ defined to be
	\begin{equation}\label{eq-G-general-np}
		H(I) = \sum_{t \in I} \mathpzc{k}(X_t, X_t) - \frac{1}{|I|} \sum_{t \in I} \sum_{s \in I}k(X_t, X_s).
	\end{equation}
	When $\mathcal{X} = \mathbb{R}$ and $\mathpzc{k}(\cdot, \cdot)$ is the linear kernel, then \eqref{eq-G-general-np} is exactly \eqref{eq-G}.

\begin{theorem}[Theorem 3.1 in \citealp{garreau2018consistent}] \label{thm-kernel-upper}
Let $\{X_t\}_{t=1}^T$ satisfy \Cref{assump-general-np}.  Assume that there exists a positive absolute constant $M$ such that 
	\begin{equation}\label{eq-k-upper-fixed}
		\mathpzc{k}(X_t, X_t) \leq M^2 < \infty, \quad \forall t \in \{1, \ldots, T\}.
	\end{equation}
    Assume that there exists a sufficiently large absolute constant $C_{\mathrm{SNR}} > 0$ such that for any $\xi > 0$,
	\[
		\kappa\sqrt{\Delta} \geq C_{\mathrm{SNR}} K\sqrt{\log^{1+\xi}(T)}.
	\]
	Let $\{\widehat{\eta}_k\}_{k = 1, \ldots, \widehat{K}}$ be the output of \Cref{algorithm:PDP} with the loss function defined in \eqref{eq-G-general-np}.  We have that, for an absolute constant $C > 0$, define $\lambda = CK^2 \log(T)$.  It holds that 
	\begin{align*}
		\mathbb{P}\Bigl\{  \widehat K =K ; \quad   |\eta_k-\hat \eta_k| \le  C_{\epsilon} K \log(T) \kappa^{-2}, \text{ for all }  k \Bigr\} \geq 1 - T^{-c}, 
	\end{align*}
	where $C_{\epsilon}, c > 0$ are absolute constants. 	
\end{theorem}

Note that the localisation error obtained in \Cref{thm-kernel-upper} is linear in $K$, the number of true change points.  In order to match the minimax lower bound $\kappa^{-2}$, we can adopt the refinement idea again.

\begin{corollary} \label{thm-refine-np-general}
Under the same settings and conditions in \Cref{thm-kernel-upper}, let $\{\nu_k\}_{k = 1}^K$ be a set of initial change point estimators satisfying $\max_{k = 1, \ldots, K}|\nu_k - \eta_k| \leq \Delta/5$.  For each $k \in \{1, \ldots, K\}$, define
	\[
		s_k = \nu_{k-1}/2 + \nu_k/2, \quad e_k = \nu_k/2 + \nu_{k+1}/2 \quad \mbox{and}\quad I_k = [s_k, e_k],
	\]
	with $\nu_0 = 1$ and $\nu_{K+1} = T + 1$.  For $k \in \{1, \ldots, K\}$, we let
	\begin{align*}
		& \widetilde{\eta}_k = \min_{t \in \{s_k + 1, \ldots, e_k - 1\}} \Bigg\{\sum_{i \in I_k} \mathpzc{k}(X_i, X_i) - \frac{1}{t-s_k} \sum_{i \in [s_k, t]} \sum_{j \in [s_k, t]} \mathpzc{k}(X_i, X_j) \\
		& \hspace{5cm} - \frac{1}{e_k-t} \sum_{i \in (t, e_k]} \sum_{j \in (t, e_k]} \mathpzc{k}(X_i, X_j)\Bigg\}.
	\end{align*}
	Then we have
	\begin{align*}
		\mathbb{P}\Bigl\{  \widehat K =K ; \quad   |\eta_k- \widetilde{\eta}_k| \le  C_{\epsilon} \log(T) \kappa^{-2}, \text{ for all }  k \Bigr\} \geq 1 - T^{-c}, 
	\end{align*}
	where $C_{\epsilon}, c > 0$ are absolute constants. 	
\end{corollary}

\Cref{thm-refine-np-general} is straightforward based on the observation that there is one and only one change point in each interval $I_k$, $k = 1, \ldots, K$.  Together with \Cref{thm-kernel-upper}, it shows that the \emph{de facto} $K$ is exactly one in each working interval $I_k$, and the results hold.

Considering \Cref{thm-refine-np-general} as a refinement step of \Cref{algorithm:PDP}, we show that a penalisation-based-method is able to achieve the near optimality in localisation error, off by a logarithmic factor.  It remains an \textbf{open question} that in terms of the signal-to-noise ratio, if one can weaken the condition from $\kappa\sqrt{\Delta} \geq C_{\mathrm{SNR}} K\sqrt{\log^{1+\xi}(T)}$ to $\kappa\sqrt{\Delta} \geq C_{\mathrm{SNR}} \sqrt{\log^{1+\xi}(T)}$.

A more interesting open problem in this problem is how one can achieve the nonparametric rate in the reproducing kernel Hilbert space change point analysis.  To explain this, we remark that despite the great flexibility in terms of $\mathcal{X}$ we reviewed in this subsection, the rates achieved are the same as the rates in \Cref{sec-detec-local}, which deals with a parametric problem.  This is because the \emph{de facto} complexity of the space is hidden in the conditions that $V$ and $M$ are absolute constants.  A more thorough result, which remains open, should involve the complexity of the reproducing kernel Hilbert space.  For instance, one should consider the Rademacher complexity \citep[e.g.][]{mendelson2002geometric, bartlett2005local} and the covering number of the reproducing kernel Hilbert space unit-ball in $\ell_2$-norm.  If the reproducing kernel Hilbert space is a Sobolev space $W^{\alpha, 2}$, $\alpha > 1/2$, we conjecture that the detection boundary should be of order $O(T^{-\alpha/(2\alpha+1)})$ and the localisation rate should be of order $O(T^{-2\alpha/(2\alpha+1)})$.  However, these remain as \textbf{open problems} in this area.

\subsection{Conclusions}\label{sec-np-conclusions}

In this section, we reviewed three different nonparametric change point detection problems.  We would like to mention that the key to deploy \Cref{algorithm:PDP} is to define a suitable loss function.  It is the sum of residual squares in Sections~\ref{sec-detec-local} and \ref{sec-piecewise-polynomial}, and a kernel version of the sum of residual squares in \Cref{sec-general-np}.  The key for \Cref{algorithm:PDP} to execute in polynomial time is that the loss function is separable in terms of the intervals.  To be specific, it should be of the form
	\begin{equation}\label{eq-G-separable-form}
		G(\mathcal{P}, \{X_t\}, \lambda) = \sum_{I \in \mathcal{P}} H(I) + \lambda |\mathcal{P}|
	\end{equation}
	and $H(I)$ is solely a function relying on the data in the interval $I$.  This suggests that it is not clear how one can directly apply \Cref{algorithm:PDP} to the problems studied in Sections~\ref{sec-np-uni} and \ref{sec-np-multi}.  For instance, in \Cref{sec-np-uni}, we see the jump is defined in \eqref{eq-as1-kappa} and it is the Kolmogorov--Smirnov distance between two different distributions.  To estimate the Kolmogorov--Smirnov distance, one needs to know where on the support the difference is taken to be largest, namely 
	\[
		z^* \in \argmax_{z \in \mathbb{R}}|F_1(z) - F_2(z)|.
	\]  
	Back to \eqref{eq-G-separable-form}, in order to estimate the change points, a certain form of the loss function is inevitably a function of $z^*$, which is not solely determined by one interval.  

The reason that \Cref{algorithm:PDP} is applicable in \Cref{sec-general-np} is due to the construction of the reproducing kernel Hilbert space, which transforms the change points of the distributions of data, to the change points in the Bochner integrals.  This transformation to a certain extent turns a general nonparametric problem to a parametric one, which echos the discussions at the end of \Cref{sec-general-np}.

Finally, we mention a list of papers working on different aspects in this area: \cite{hawkins2010nonparametric}, \cite{haynes2017computationally-2}, \cite{itoh2010change}, \cite{matteson2014nonparametric}, \cite{vanegas2019multiscale}, \cite{zou2014nonparametric}, \cite{harchaoui2007retrospective}, \cite{garreau2018consistent}, \cite{celisse2018new} and \cite{arlot2019kernel}, among others.

\section{Conclusions}

In this survey, we covered a range of change point analysis problems, focusing on the minimax rates of detection and localisation, with an emphasis on distinguishing these two ideas.  The univariate mean change point detection problem lays down the foundation in terms of minimax lower bounds and two types of popular methods which are nearly optimal.  For more complicated cases we covered in this survey, we have reviewed different situations, in some cases we can show the near optimality in terms of both detection and localisation, in some cases we show that an extra refinement step can reach the near optimality in localisation but the detection conditions remain sub-optimal, and in some cases we show that the near optimality is reachable under some extra conditions.

There are still many open questions in the change point analysis area.  Throughout the survey, we have identified a few.  In addition, we would like to reiterate that the minimax results we reviewed in this paper are all on offline change point detection and localisation.  The minimax rates of online change point detection and localisation, and minimax rates of both online and offline change point testing remain largely unknown.  The optimality we achieved in this survey are all off by logarithmic factors.  It would also be interesting to further refine the results improving the results.  Efforts along this line include some results in \cite{verzelen2020optimal}.

\bibliographystyle{ims}
\bibliography{ref}

\end{document}